\documentclass[12pt]{amsart}

\usepackage{geometry}
\usepackage{amsfonts, adjustbox, amssymb, amscd}
\numberwithin{equation}{section}

\usepackage{bm}
\usepackage{verbatim}
\usepackage{mathrsfs}
\usepackage{graphicx}
\usepackage{tikz-cd}
\usepackage{subcaption}
\usepackage{listings}
\usepackage{subfiles}
\usepackage[toc,page]{appendix}
\usepackage{mathtools}
\usepackage{comment}
\usepackage{enumerate}
\usepackage{enumitem}
\usepackage[all]{xy}

\usepackage{graphicx}
\usepackage{appendix}
\usepackage{hyperref}
\hypersetup{
    colorlinks=true,
    citecolor=red,
    linkcolor=blue,
    filecolor=magenta,      
    urlcolor=red,
}
\newcommand{\bb}{\bm{b}}

\newcommand{\Mm}{{\bf{M}}}
\newcommand{\Nn}{{\bf{N}}}

\newcommand{\NN}{{\bf{N}}}
\newcommand{\Dd}{{\bf{D}}}

\newcommand{\Qq}{\mathbb{Q}}

\newcommand{\Rr}{\mathbb{R}}
\newcommand{\RR}{{\bf{R}}}

\newcommand{\Span}{\operatorname{Span}}

\newcommand{\Center}{\operatorname{center}}

\newcommand{\Exc}{\operatorname{Exc}}

\newcommand{\rk}{\operatorname{rank}}

\newcommand{\ninv}{\operatorname{ninv}}
\newcommand{\inv}{\operatorname{inv}}

\newcommand{\tmld}{{\operatorname{tmld}}}

\newcommand{\Weil}{\operatorname{Weil}}

\newcommand{\Supp}{\operatorname{Supp}}

\newcommand{\mult}{\operatorname{mult}}

\newcommand{\Aa}{{\mathfrak{A}}}

\newcommand{\Ff}{\mathcal{F}}

\newcommand{\MM}{{\mathfrak{M}}}

\newcommand{\Ii}{\Gamma}

\newcommand{\Ll}{{\bf{L}}}


\newcounter{parentnumber}

\newtheorem{thm}{Theorem}[section]

\newtheorem{lem}[thm]{Lemma}

\theoremstyle{definition}
\newtheorem{defn}[thm]{Definition}

\theoremstyle{definition}

\newtheorem{defthm}[thm]{Definition-Theorem}
\newtheorem{setup}[thm]{Set-up}
\newtheorem{Asp}[thm]{Assumption}

\newtheorem{nota}[thm]{Notation}

\newtheorem{cons}[thm]{Construction}

\theoremstyle{definition}

\begin{document}

\title{Sarkisov program for algebraically integrable and threefold foliations}
\author{Yifei Chen}
\author{Jihao Liu}
\author{Yanze Wang}

\date{\today}

\begin{abstract}
By applying the theory of the minimal model program for adjoint foliated structures, we establish the Sarkisov program for algebraically integrable foliations on klt varieties: any two Mori fiber spaces of such structure are connected by a sequence of Sarkisov links. Combining with a result of R. Mascharak, we establish the Sarkisov program for foliations in dimension at most $3$ with mild singularities. Log version and adjoint foliated version of the aformentioned Sarkisov programs are also established.
\end{abstract}

\address{State Key Laboratory of Mathematical Sciences, Academy of Mathematics and Systems Science, Chinese Academy of Sciences, Beijing 100190, China}
\email{yifeichen@amss.ac.cn}

\address{Department of Mathematics, Peking University, No. 5 Yiheyuan Road, Haidian District, Peking 100871, China}
\email{liujihao@math.pku.edu.cn}

\address{Academy of Mathematics and Systems Science, Chinese Academy of Sciences, No. 55 Zhonguancun East Road, Haidian District, Beijing, 100190, China}
\email{wangyanze@amss.ac.cn}

\maketitle

\pagestyle{myheadings}\markboth{\hfill Yifei Chen, Jihao Liu, and Yanze Wang\hfill}{\hfill Sarkisov program for algebraically integrable and threefold foliations\hfill}

\tableofcontents
%
%

\section{Introduction}\label{sec:Introduction}
We work over the field of complex numbers $\mathbb{C}$.

\medskip

\noindent\textbf{Sarkisov program of varieties.} Let $X$ be a uniruled smooth projective variety. The minimal model program (MMP) predicts that, after running a $K_X$-MMP, that is, a sequence of $K_X$-divisorial contractions and flips $X \dashrightarrow X'$, one should obtain a Mori fiber space $X' \to Z$, i.e. a fibration such that $\dim X' > \dim Z$, $-K_{X'}$ is ample$/Z$, and $\rho(X'/Z) = 1$. The existence of such a fibration follows from running a $K_X$-MMP with scaling of an ample divisor (cf. \cite[Corollary 1.4.2]{BCHM10}), though the resulting Mori fiber space is generally not unique. For instance, the Hirzebruch surface $X := \mathbb{F}_1$ admits a natural $\mathbb{P}^1$-fibration structure which is a $K_X$-Mori fiber space, yet there also exists a $K_X$-MMP step $X \to \mathbb{P}^2$, where $\mathbb{P}^2 \to \{\mathrm{pt}\}$ is another Mori fiber space.

This non-uniqueness naturally raises the question of how different Mori fiber spaces of $X$ are related. The study of such connections, now known as the \emph{Sarkisov program}, was initiated by Sarkisov \cite{Sar80,Sar82} in his work on conic bundles over threefolds. The theory was further developed by Corti \cite{Cor95}, Iskovskikh \cite{Isk96}, and Bruno-Matsuki \cite{BM97}, among others, in the context of the MMP for klt pairs. Roughly speaking, the Sarkisov program aims to decompose any birational map $X_1 \dashrightarrow X_2$ between Mori fiber spaces into a sequence of elementary birational transformations called \emph{Sarkisov links} (see Definition \ref{defn: sarkisov links} for a precise description). This decomposition plays an important role in the study of the birational automorphism groups of Mori fiber spaces associated with Fano varieties and the study of higher-dimensional Cremona groups (cf. \cite{Isk91,IKT93,Isk96,LZ20,BLZ21}).

Following the proof of the existence of the minimal model program for klt varieties \cite{BCHM10}, Hacon and M\textsuperscript{c}Kernan \cite[Theorem 1.1]{HM13} showed that any two Mori fiber spaces $X_1 \to Z_1$ and $X_2 \to Z_2$ that are outputs of a $K_X$-MMP are connected by a sequence of Sarkisov links. An alternative proof was given in an unpublished note of Hacon \cite{Hac12} which was later reformulated to obtain the Sarkisov program generalized pairs by the second author \cite{Liu21} (see \cite{CW25} by the first and the third author for a survey). Variations of the Sarkisov program have since been developed in several directions, such as the group-equivariant Sarkisov program \cite{Flo20} and the Sarkisov program in positive and mixed characteristics for low-dimensional varieties \cite{BFSZ24,Sti21}. For further related work, we refer the reader to \cite{Kal13,Miy19,FP22,He24}.

\medskip

\noindent\textbf{Sarkisov program for foliations.} In recent years, there are substantial progress in establishing the foundations of the minimal model program (MMP) for foliations. For an lc foliation $\mathcal{F}$ of rank $r$ on a klt variety $X$ of dimension $d$, the existence of a $K_{\mathcal{F}}$-MMP (i.e. the cone theorem, contraction theorem, and existence of flips) has been established when $d=3$ \cite{McQ08, Bru15}, when  $(r,d)=(2,3)$ and $\mathcal{F}$ is $\mathbb{Q}$-factorial F-dlt \cite{Spi20, CS21, SS22}, when $(r,d)=(1,3)$ and $X$ is $\mathbb{Q}$-factorial \cite{CS20}, and when $\mathcal{F}$ is algebraically integrable, regardless of $(r,d)$ \cite{CHLX23, LMX24, CHLMSSX24}.

Given these developments, it becomes natural to investigate whether the Sarkisov program admits a counterpart for foliations. Just as the classical Sarkisov program is deeply connected to the study of birational automorphism groups of Fano varieties, its foliated version is expected to play an important role in understanding birational automorphism groups of foliations. A recent work by R. Mascharak \cite[Theorem 1.1]{Mas24} has established the Sarkisov program for foliations for the case $(r,d)=(2,3)$, assuming that $\Ff$ has mild singularities ($\mathbb Q$-factorial foliated F-dlt). While not explicitly stated, the methods employed in the proof extend naturally to the $(r,d)=(1,2)$ case. However, the $(r,d)=(1,3)$ case remains largely open due to the failure of Bertini-type theorems, despite some partial results, e.g. \cite[Theorem 1.3]{Mas24}.

The first goal of this paper is to establish the Sarkisov program for rank $1$ lc foliations on $\mathbb Q$-factorial klt threefolds. 

\begin{thm}\label{thm: main dim 3 rank 1}
Sarkisov program holds for lc foliations of rank $1$ on $\mathbb Q$-factorial projective klt threefolds. More precisely, let $\Ff$ be an lc foliation of rank $1$ on a $\mathbb Q$-factorial projective klt threefold $X$. For $i\in\{1,2\}$, let $(X,\Ff)\dashrightarrow (X_i,\Ff_i)$ be two $K_{\Ff}$-MMPs and let $f_i: X_i\rightarrow Z_i$ be $K_{\Ff_i}$-Mori fiber spaces. Then the $f_1$ and $f_2$ are connected by a finite sequence of Sarkisov links.
\end{thm}

Combining with \cite[Theorem 1.1]{Mas24}, we obtain the Sarkisov program in dimension $\leq 3$ for foliations with mild singularities. Here ``mild" stands for ``lc foliations on $\mathbb Q$-factorial klt varieties" unless $\rk\Ff=2$ and $\dim X=3$, and stands for ``$\mathbb Q$-factorial F-dlt foliations" (\cite[Definition 3.1]{CS21}) if $\rk\Ff=2$ and $\dim X=3$. In either case, these foliations are always outputs of foliated log smooth foliations. 

We have the following precise result:

\begin{thm}\label{thm: main dim 3}
Sarkisov program holds for foliated log smooth foliations in dimension $\leq 3$.
\end{thm}

Theorem \ref{thm: main dim 3 rank 1} is actually a very special case of the following result on the Sarkisov program for $\mathbb Q$-factorial klt algebraically integrable adjoint foliated structures. 

\begin{thm}\label{thm: main afs}
Sarkisov program holds for lc algebraically integrable adjoint foliated structures (see Definition \ref{defn: afs}) on $\mathbb Q$-factorial klt varieties. More precisely: 

Let $\Aa_W/U$ be an lc algebraically integrable adjoint foliated structure whose ambient variety $W$ is $\mathbb Q$-factorial klt. Let $g: W\dashrightarrow X$, $g': W\dashrightarrow X'$ be two $K_{\Aa}$-MMPs$/U$, $\Aa:=g_*\Aa_W$, $\Aa':=g'_*\Aa_W$, and let $f: X\rightarrow Z$, $f': X'\rightarrow Z$ be a $K_{\Aa}$-Mori fiber space$/U$ and a  $K_{\Aa'}$-Mori fiber space$/U$ respectively. Then the $f$ and $f'$ are connected by a finite sequence of Sarkisov links$/U$.
\end{thm}

As a special case of Theorem \ref{thm: main afs}, we obtain:
\begin{thm}\label{thm: main algint}
    Sarkisov program holds for lc algebraically integrable foliations on $\mathbb Q$-factorial klt varieties.
\end{thm}

We outline the key ideas in the proof of Theorem \ref{thm: main dim 3 rank 1}. The main challenge in establishing the Sarkisov program for rank $1$ foliations on threefolds is the failure of Bertini-type theorems \cite{Mas24}. While recent advances in the minimal model program for foliations suggest that such failures can often be circumvented by working with generalized foliated quadruples \cite{LLM23,CHLX23,LMX24} via the simple trick of incorporating the auxiliary ample divisor into the nef part, this approach is ineffective for the Sarkisov program. There are two major difficulties: First, the foliation in the Sarkisov program cannot be assumed to be klt, so the classical strategy of \cite{HM13} is inapplicable and we do not have a good understanding on the geometry of ample models of foliations. Second, the alternative approach ``Sarkisov program with double scaling" approach from \cite{Liu21} fails because terminalization usually do not exist for lc foliations.

To overcome these difficulties, in this paper, we shall use the theory of \emph{adjoint foliated structures} instead. The crucial observation is that any $K_{\mathcal{F}}$-MMP can be viewed as a $(tK_{\mathcal{F}} + (1-t)K_X + A)$-MMP for some $0 < 1-t \ll 1$ and an ample $\mathbb{R}$-divisor $A$, where $(X, \mathcal{F}, 0, \overline{A}, t)$ forms a klt adjoint foliated structure. In this case, we can replace the category of ``lc foliations on klt varieties" with the category of ``klt adjoint foliated structures". The later is a more technical and complicated category, but thanks to the recent development on the minimal model program for algebraically integrable adjoint foliated structures \cite{CHLMSSX24,CHLMSSX25}, we can apply the arguments as in \cite{Liu21} to deduce the corresponding Sarkisov program (Theorem \ref{thm: main afs}). We conclude by noting that all rank $1$ foliations appearing in the Sarkisov program for potentially klt varieties have non-pseudo-effective canonical bundles and are consequently algebraically integrable \cite{CP19,ACSS21}.

\medskip

\noindent\textbf{Acknowledgement.} The work is supported by the National Key R\&D Program of China (\#2024YFA1014400). The first and third authors are partially supported by the NSFC grant No. 12271384. The authors would like to thank Roktim Mascharak for useful discussions.

\section{Preliminaries}

We will adopt the standard notation and definitions on MMP in \cite{KM98,BCHM10} and use them freely. For adjoint foliated structures, generalized foliated quadruples, foliated triples, foliations, and generalized pairs, we adopt the notation and definitions in \cite{CHLMSSX24,CHLMSSX25}  which generally align with \cite{LLM23,CHLX23} (for generalized foliated quadruples), \cite{CS20, ACSS21, CS21} (for foliations and foliated triples), and \cite{BZ16,HL23} (for generalized pairs and $\bb$-divisors), possibly with minor differences. 

\subsection{Notation}

\begin{defn}
    A \emph{contraction} is a projective morphism of varieties
    $f \colon X \to Y$
    such that 
    $f_\ast \mathcal{O}_X=\mathcal{O}_Y$. When $X, Y$ are normal, the condition 
$f_\ast \mathcal{O}_X=\mathcal{O}_Y$
is equivalent to the fibers of $f$ being connected.
\end{defn}

\begin{nota}
    Let $h \colon  X\dashrightarrow X'$ be a birational map between normal varieties. By abusing notation, we denote by $\Exc(h)$ the reduced divisor supported on the codimension one part of the exceptional locus of $h$. 
\end{nota}

\begin{defn}[NQC]
    Let $X\rightarrow U$ be a projective morphism between normal quasi-projective varieties. Let $D$ be a nef $\Rr$-Cartier $\Rr$-divisor on $X$ and $\Mm$ a nef $\bb$-divisor on $X$. 
        We say that $D$ is \emph{NQC}$/U$ if $D=\sum d_iD_i$, where each $d_i\geq 0$ and each $D_i$ is a nef$/U$ Cartier divisor. Here NQC stands for ``Nef $\mathbb Q$-Cartier Combinations", cf. \cite[Definition 2.15]{HL22}. We say that $\Mm$ is \emph{NQC}$/U$ if $\Mm=\sum \mu_i\Mm_i$, where each $\mu_i\geq 0$ and each $\Mm_i$ is a nef$/U$ $\bb$-Cartier $\bb$-divisor.
\end{defn}

\subsection{Adjoint foliated structures}

\begin{defn}[Foliations, {cf. \cite{ACSS21,CS21}}]\label{defn: foliation}
Let $X$ be a normal variety. A \emph{foliation} on $X$ is a coherent sheaf $\Ff\subset T_X$ such that
\begin{enumerate}
    \item $\Ff$ is saturated in $T_X$, i.e. $T_X/\Ff$ is torsion free, and
    \item $\Ff$ is closed under the Lie bracket.
\end{enumerate}

The \emph{rank} of a foliation $\Ff$ on a variety $X$ is the rank of $\Ff$ as a sheaf and is denoted by $\rk\Ff$. 
The \emph{co-rank} of $\Ff$ is $\dim X-\rk\Ff$. The \emph{canonical divisor} of $\Ff$ is a divisor $K_\Ff$ such that $\mathcal{O}_X(-K_{\mathcal{F}})\cong\mathrm{det}(\Ff)$. If $\Ff=0$, then we say that $\Ff$ is a \emph{foliation by points}.

Given any dominant map 
$h: Y\dashrightarrow X$ and a foliation $\mathcal F$ on $X$, we denote by $h^{-1}\Ff$ the \emph{pullback} of $\Ff$ on $Y$ as constructed in \cite[3.2]{Dru21} and say that $h^{-1}\Ff$ is \emph{induced by} $\Ff$. Given any birational map $g: X\dashrightarrow X'$, we denote by $g_\ast \Ff:=(g^{-1})^{-1}\Ff$ the \emph{pushforward} of $\Ff$ on $X'$ and also say that $g_\ast \Ff$ is \emph{induced by} $\Ff$. We say that $\Ff$ is an \emph{algebraically integrable foliation} if there exists a dominant map $f: X\dashrightarrow Z$ such that $\Ff=f^{-1}\Ff_Z$, where $\Ff_Z$ is the foliation by points on $Z$, and we say that $\Ff$ is \emph{induced by} $f$.

A subvariety $S\subset X$ is called \emph{$\Ff$-invariant} if for any open subset $U\subset X$ and any section $\partial\in H^0(U,\Ff)$, we have $\partial(\mathcal{I}_{S\cap U})\subset \mathcal{I}_{S\cap U}$,  where $\mathcal{I}_{S\cap U}$ denotes the ideal sheaf of $S\cap U$ in $U$.  
For any prime divisor $P$ on $X$, we define $\epsilon_{\Ff}(P):=1$ if $P$ is not $\Ff$-invariant and $\epsilon_{\Ff}(P):=0$ if $P$ is $\Ff$-invariant. For any prime divisor $E$ over $X$, we define $\epsilon_{\Ff}(E):=\epsilon_{\Ff_Y}(E)$ where $h: Y\dashrightarrow X$ is a birational map such that $E$ is on $Y$ and $\Ff_Y:=h^{-1}\Ff$.

If foliation structure $\Ff$ on $X$ is clear in the context, then given an $\mathbb R$-divisor $D = \sum_{i = 1}^k a_iD_i$ where each $D_i$ is a prime divisor,
we denote by $D^{{\rm ninv}} \coloneqq \sum \epsilon_{\mathcal F}(D_i)a_iD_i$ the \emph{non-$\Ff$-invariant part} of $D$ and $D^{{\rm inv}} \coloneqq D-D^{{\rm ninv}}$ the \emph{$\Ff$-invariant part} of $D$.
\end{defn}

\begin{defn}\label{defn: afs}
An \emph{adjoint foliated structure} $\Aa/U:=(X,\Ff,B,\Mm,t)/U$ is the datum of a normal quasi-projective variety $X$, a projective morphism $X\rightarrow U$, a foliation $\Ff$ on $X$, an $\Rr$-divisor $B\geq 0$ on $X$, a $\bb$-divisor $\Mm$ nef$/U$, and a real number $t\in [0,1]$ such that $$K_{\Aa}:=K_{(X,\Ff,B,\Mm,t)}:=tK_{\Ff}+(1-t)K_X+B+\Mm_X$$ is $\Rr$-Cartier. 
$K_{\Aa}$ is called as the \emph{canonical $\Rr$-divisor} of $\Aa$. $X,\Ff,B,\Mm,t$ are called the \emph{ambient variety}, \emph{foliation part}, \emph{boundary part}, \emph{nef part} (or \emph{moduli part}), and \emph{parameter} of $\Aa$ respectively. The $\mathbb R$-divisor $B+\Mm_X$ is called the \emph{generalized boundary} of $\Aa$. We say that $\Aa/U$ is \emph{of general type} if $K_{\Aa}$ is big$/U$. 

For any $\Rr$-divisor $D$ on $X$ and $\bb$-divisor $\Nn$ on $X$ such that $D+\Nn_X$ is $\Rr$-Cartier and $\Mm+\Nn$ is nef$/U$, we denote by $(\Aa,D,\Nn):=(X,\Ff,B+D,\Mm+\Nn,t)$. If $D=0$ then we may drop $D$, and if $\Nn=\bm{0}$ then we may drop $\Nn$. For any projective birational morphism $h: X'\rightarrow X$, we define 
$$h^*\Aa:=(X',\Ff',B',\Mm,t)$$
where $\Ff':=h^{-1}\Ff$ and $B'$ is the unique $\Rr$-divisor such that $K_{h^*\Aa}=h^*K_{\Aa}$. 
For any birational map$/U$ $\phi: X\dashrightarrow X'$ which does not extract any divisor and any adjoint foliated structure $\Aa/U$ on $X$, we define $\phi_*\Aa/U$ to be the adjoint foliated structure such that $\phi_*\Aa:=(X',\phi_*\Ff,\phi_*B,\Mm,t)$ and say that $\phi_*\Aa$ is the \emph{image} of $\Aa$ on $X'$. For any prime divisor $E$ on $X'$, we denote by
$$a(E,\Aa):=-\mult_EB'$$
the \emph{discrepancy} of $E$ with respect to $\Aa$. The \emph{total minimal log discrepancy} of $\Aa$ is
$$\tmld(\Aa):=\inf\{a(E,\Aa)+t\epsilon_{\Ff}(E)+(1-t)\mid E\text{ is over }X\}$$
We say that $\Aa$ is \emph{lc} (resp. \emph{klt}) if $\tmld(\Aa)\geq 0$ (resp. $>0$). We say that $\Aa$ is \emph{terminal} if $a(E,\Aa)>0$ for any prime divisor $E$ that is exceptional$/X$. If we allow $B$ to have negative coefficients, then we add the prefix ``sub-" for the types of singularities above.

When $\Mm=\bm{0}$, and either $t=0$ or $\Ff=T_X$, we call $(X,B)/U$ a \emph{pair}. We say that $(X,\Ff,B,\Mm,t)/U$ is \emph{NQC} if $\Mm$ is NQC$/U$. If $B=0$, or if $\Mm=\bm{0}$, or if $U$ is not important, then we may drop $B,\Mm,U$ respectively. If $U=\{pt\}$ then we also drop $U$ and say that $(X,\Ff,B,\Mm,t)$ is \emph{projective}. If we allow $B$ to have negative coefficients, then we shall add the prefix ``sub-". If $B$ is a $\Qq$-divisor, $\Mm$ is a $\Qq$-$\bb$-divisor, and $t\in\mathbb Q$, then we shall add the prefix ``$\Qq$-".
\end{defn}

\subsection{Birational maps in the MMP}

\begin{defn}
    Let $X\rightarrow U$ be a projective morphism from a normal quasi-projective variety to a variety.  Let $D$ be an $\Rr$-Cartier $\Rr$-divisor on $X$ and $\phi: X\dashrightarrow X'$ a birational map$/U$. Then we say that $X'$ is a \emph{birational model} of $X$. We say that $\phi$ is $D$-non-positive (resp. $D$-negative, $D$-trivial, $D$-non-negative, $D$-positive) if the following conditions hold:
    \begin{enumerate}
    \item $\phi$ does not extract any divisor.
    \item $D':=\phi_\ast D$ is $\Rr$-Cartier.
    \item There exists a resolution of indeterminacy $p: W\rightarrow X$ and $q: W\rightarrow X'$, such that
    $$p^\ast D=q^\ast D'+F$$
    where $F\geq 0$ (resp. $F\geq 0$ and $\Supp p_\ast F=\Exc(\phi)$, $F=0$, $0\geq F$, $0\geq F$ and $\Supp p_\ast F=\Exc(\phi)$).
    \end{enumerate}
\end{defn}

\begin{defn}
    Let  $f: X\rightarrow Z$ be a contraction$/U$ between normal quasi-projective varieties, and $D$ an $\mathbb R$-Cartier $\mathbb R$-divisor on $X$. $f$ is called a \emph{$D$-Mori fiber space} if $\dim X>\dim Z$, $\rho(X/Z)=1$, and $D$ is anti-ample$/Z$. If $f$ is a $D$-Mori fiber space for some $D$, then we say that $f$ is a \emph{Mori fiber space}.
\end{defn}

\begin{defn} Let $X\rightarrow U$ be a projective morphism from a normal quasi-projective variety to a variety. Let $D$ be an $\Rr$-Cartier $\Rr$-divisor on $X$, $\phi: X\dashrightarrow X'$ a $D$-non-positive birational map$/U$, and $D':=\phi_\ast D$.
\begin{enumerate}
    \item If $D'$ is nef$/U$, then $X'$ is called a \emph{weak lc model}$/U$ of $D$.
    \item If $D'$ is nef$/U$ and $\phi$ is $D$-negative, then $X'$ is called a \emph{minimal model}$/U$ of $D$.
    \item If $\phi$ is $D$-negative and $D'$ is semi-ample$/U$, then $X'$ is called a \emph{good minimal model}$/U$ of $D$.
    \item Assume that there exists a contraction$/Z$ $f: X'\rightarrow Z$ that is a $D'$-Mori fiber space, and $\phi$ is $D$-negative. Then $f$ is called a \emph{Mori fiber space}$/U$ of $D$. 
\end{enumerate}
\end{defn}

\begin{defn}\label{defn: minimal model}
Let $\Aa/U$ be an adjoint foliated structure with ambient variety $X$ and let $\phi: X\dashrightarrow X'$ be a birational map$/U$. Let $\Aa':=\phi_*\Aa$. We say that $\Aa'/U$ is weak lc model (resp. minimal model, good minimal model) of $\Aa/U$ if $X'$ is a weak lc model$/U$ (resp. minimal model$/U$, good minimal model$/U$) of $K_{\Aa}$. A contraction$/U$ $f: X'\rightarrow Z$ is called a \emph{Mori fiber space} of $\Aa/U$ if $f$ is a Mori fiber space$/U$ of $K_{\Aa}$. We say that $f$ is $\mathbb Q$-factorial if $X'$ is $\mathbb Q$-factorial.
\end{defn}

\begin{defn}[{cf. \cite[3.2 Log canonical foliated pairs]{ACSS21}, \cite[Definition 6.2.1]{CHLX23}}]\label{defn: foliated log smooth}
Let $\Aa/U:=(X,\Ff,B,\Mm,t)/U$ be an algebraically integrable adjoint foliated structure. We say that $\Aa$ is \emph{foliated log smooth} if there exists a contraction $f: X\rightarrow Z$ satisfying the following.
\begin{enumerate}
  \item $X$ has at most quotient toric singularities.
  \item $\Ff$ is induced by $f$.
  \item $(X,\Sigma_X)$ is toroidal for some reduced divisor $\Sigma_X$ such that $\Supp B\subset\Sigma_X$.  In particular, $(X,\Supp B)$ is toroidal, and $X$ is $\Qq$-factorial klt.
  \item There exists a log smooth pair $(Z,\Sigma_Z)$ such that $$f: (X,\Sigma_X)\rightarrow (Z,\Sigma_Z)$$ is an equidimensional toroidal contraction.
  \item $\Mm$ descends to $X$.
\end{enumerate}
We say that $f: (X,\Sigma_X)\rightarrow (Z,\Sigma_Z)$ is \emph{associated with} $(X,\Ff,B,\Mm,t)$. It is important to remark that $f$ may not be a contraction$/U$. In particular, $\Mm$ may not be nef$/Z$.

Note that the definition of foliated log smooth has nothing to do with $t$. In other words, as long as $(X,\Ff,B,\Mm)$ is foliated log smooth, $(X,\Ff,B,\Mm,t)$ is foliated log smooth.
\end{defn}

\begin{defn}[Foliated log resolutions]\label{defn: log resolution}
Let $\Aa/U:=(X,\Ff,B,\Mm,t)/U$ be an algebraically integrable sub-adjoint foliated structure. A \emph{foliated log resolution} of $\Aa$ is a birational morphism $h: X'\rightarrow X$ such that 
$$(X',\Ff':=h^{-1}\Ff,B':=h^{-1}_\ast B+\Exc(h),\Mm,t)$$ 
is foliated log smooth. By \cite[Lemma 6.2.4]{CHLX23}, foliated log resolution for $\Aa$ always exists.
\end{defn}

\subsection{Terminalization}

\begin{lem}\label{lem: klt t<1}
    Let $\Aa/U=(X,\Ff,B,\Mm,t)$ be an algebraically integrable sub-adjoint foliated structure. If $t=1$ and $\Ff\not=T_X$, then $\Aa$ is not sub-klt.
\end{lem}
\begin{proof}
There exists an $\Ff$-invariant divisor $D$ on $X$ such that $D$ is not a component of $\Supp B$. We have $a(D,\Aa)=0$, so $\Aa$ is not sub-klt.
\end{proof}

\begin{lem}\label{lem: finite non-terminal place}
    Let $\Aa/U$ be a sub-klt algebraically integrable sub-adjoint foliated structure with ambient variety $X$. Then
    $$\mathcal{S}_0:=\{F\mid F\text{ is exceptional}/X,a(F,\Aa)\leq 0\}$$
    is a finite set.
\end{lem}
\begin{proof}
  Possibly replacing $\Aa$ with a foliated log resolution, we may assume that $\Aa=(X,\Ff,B,\Mm,t)$ is foliated log smooth, and then we may assume that $\Mm=\bm{0}$. We may assume that $\dim X\geq 2$. 
  
  Since $\Aa$ is klt, we may assume that $t<1$, otherwise, by Lemma \ref{lem: klt t<1}, $\Ff=T_X$ and we are done. Let $f: (X,\Sigma_X)\rightarrow (Z,\Sigma_Z)$ be a contraction associated with $(X,\Ff,B,t)$ and let $\Sigma_X^h,\Sigma_X^v$ be the horizontal$/Z$ and vertical$/Z$ part respectively. Since $\Aa$ is klt, there exists a real number $\epsilon\in (0,1)$ such that 
 $$(1-\epsilon)(\Sigma_X^h+(1-t)\Sigma_X^v)\geq B.$$
 Possibly replacing $B$, we may assume that $(1-\epsilon)(\Sigma_X^h+(1-t)\Sigma_X^v)=B$. Then we have
 $$K_{\Aa}=t(K_{\Ff}+(1-\epsilon)\Sigma_X^h)+(1-t)(K_X+(1-\epsilon)\Sigma_X).$$
 Let $$\Delta:=(1-\epsilon)\Sigma_X^h+(1-\epsilon+t\epsilon)\Sigma_X^v,$$
 then $(X,\Delta)$ is klt, and by \cite[Proposition 3.6]{ACSS21}, we have
 $$K_X+\Delta\sim_{\mathbb R,Z}K_{\Aa}.$$
Note that the set
$$\mathcal{S}:=\{E\mid E\text{ is exceptional}/X,a(E,X,\Delta)\leq 0\}$$
    is a finite set.
 Now for any prime divisor $F$ that is exceptional$/X$ and $a(F,\Aa)\leq 0$, there are two possibilities:

\medskip

\noindent\textbf{Case 1.} $V:=\Center_XF$ dominates $Z$. Let $g: X'\rightarrow X$ be a birational morphism such that $F$ is on $X'$, $\Aa':=g^*\Aa$, and $K_{X'}+\Delta':=g^*(K_X+\Delta)$. Then over the generic point of $V$ we have $K_{X'}+\Delta'=K_{\Aa'}$ and $K_X+\Delta=K_{\Aa}$,  hence
$$a(F,\Aa)=a(F,\Aa')=a(F,X',\Delta')=a(F,X,\Delta).$$
Therefore, $F\in\mathcal{S}$.

\medskip

\noindent\textbf{Case 2.}  $V:=\Center_XF$ does not dominate $Z$. Then $F$ is $\Ff$-invariant. Since $(X,\Ff,B)$ is foliated log smooth, $a(F,X,\Ff,(1-\epsilon)\Sigma_X^h)\geq 0$. Since $a(F,\Aa)\leq 0$, 
$$a(F,X,\Delta)\leq a(F,X,(1-\epsilon)\Sigma_X)\leq 0$$ 
by linearity of discrepancies. Thus $F\in\mathcal{S}$.

\medskip

So $\mathcal{S}_0\subset\mathcal{S}$ is a finite set and we are done.
\end{proof}

\begin{lem}\label{lem: perturbation terminal}
    Let $\Aa/U$ be a terminal and klt algebraically integrable sub-adjoint foliated structure with ambient variety $X$ and parameter $t<1$, let and $D\geq 0$ an $\mathbb R$-Cartier $\mathbb R$-divisor on $X$. Then $(\Aa,\epsilon D)$ is terminal and klt for any $0<\epsilon\ll 1$.
\end{lem}
\begin{proof}
    By \cite[Lemma 3.25]{CHLMSSX25}, there exists $\epsilon_0>0$ such that $(\Aa,\epsilon_0 D)$ is klt. By Lemma \ref{lem: finite non-terminal place},     
    $$\mathcal{S}_0:=\{F\mid F\text{ is exceptional}/X,a(F,\Aa,\epsilon_0D)\leq 0\}$$
    is a finite set. For any $F\in\mathcal{S}_0$, there exists $\epsilon_F>0$ such that $a(F,\Aa,\epsilon_F D)>0$. Now $(\Aa,\epsilon D)$ is terminal for any $0<\epsilon<\min_F\{\epsilon_F\}$.
\end{proof}

\begin{defthm}\label{defthm: terminalization}
  Let $\Aa/U$ be a klt algebraically integrable sub-adjoint foliated structure with ambient variety $X$. Then there exists a projective birational morphism $h: X'\rightarrow X$ such that $\Aa':=h^*\Aa$ is $\mathbb Q$-factorial terminal. $h: \Aa'\rightarrow\Aa$ is called a \emph{$\mathbb Q$-factorial terminalization} of $\Aa$. We also say that $\Aa'$ is $\mathbb Q$-factorial terminalization of $\Aa$.
\end{defthm}
\begin{proof}
    By Lemma \ref{lem: finite non-terminal place}, $\mathcal{S}_0:=\{F\mid F\text{ is exceptional}/X,a(F,\Aa)\leq 0\}$
    is a finite set. By \cite[Theorem 2.2.3]{CHLMSSX25}, there exists a projective birational morphism $h: X'\rightarrow X$ such that $X'$ is $\Qq$-factorial and the divisors contracted by $h$ are exactly the divisors in $\mathcal{S}_0$. Therefore $h$ satisfies our requirements. 
\end{proof}
\begin{lem}\label{lem: liu21 3.3}
Let $\Aa/U$ and $\Aa'/U$ be two algebraically integrable adjoint foliated structures with ambient varieties $X$ and $X'$ respectively associated with birational map $\phi: X\dashrightarrow X'$. Let $h: \Aa''\rightarrow \Aa'$ be a $\mathbb Q$-factorial terminalization of $\Aa'$ and $X''$ the ambient variety of $\Aa''$.

Assume that $K_{\Aa'}$ is nef$/U$, $\phi$ does not extract any divisor, $\Aa$ is terminal, and $a(D,\Aa)\leq a(D,\Aa')$ for any prime divisor $E$ on $X$. Then $\Aa\geq\Aa'$, $\Aa'$ is klt, and  $\psi: X\dashrightarrow X''$ does not extract any divisor. 
\end{lem}
\begin{proof}
    Let $p: W\rightarrow X$ and $q: W\rightarrow X'$ be a resolution of indeterminacy of $\phi$ and write $p^*K_{\Aa}=q^*K_{\Aa'}+E$ for some $E\geq 0$. Since $a(D,\Aa)\leq a(D,\Aa')$ for any prime divisor $E$ on $X$, $p_*E\geq 0$. Since $E$ is anti-nef$/X$, by the negativity lemma, $E\geq 0$. Thus  $\Aa\geq\Aa'$ and $\Aa'$ is klt. 

    Suppose that $\psi$ extracts a prime divisor $D$, then $D$ is exceptional$/X$ and exceptional$/X'$. Since $\Aa$ is terminal,
    $$0<a(D,\Aa)\leq a(D,\Aa')\leq 0$$
    which is not possible.
\end{proof}

\begin{lem}\label{lem: lift mfs mmp}
    Let  $\Aa/U$ be a klt algebraically integrable adjoint foliated structure and let $f_i: X_i\rightarrow Z_i$, $i\in\{1,2\}$ be two $\mathbb Q$-factorial Mori fiber spaces of $\Aa/U$ associated with birational maps $\phi_i: X\dashrightarrow X_i$. Let $\Aa_i:=(\phi_i)_*\Aa$.
    
    Then there exists a $\mathbb Q$-factorial terminal algebraically integrable adjoint foliated structure $\bar\Aa/U$ with ambient variety $\bar X$ and projective birational morphisms $h_i: \overline{X}\rightarrow X_i$, $i\in\{1,2\}$, such that each $ M $ descends to $W$, and $h_i$ can be decomposed into a sequence of steps of a $K_{\bar\Aa}$-MMP$/U$, and $(h_i)_*\bar\Aa=\Aa_i$. In particular, each $h_i$ is $K_{\bar\Aa}$-negative.
\end{lem}
\begin{proof}
    Write $\Aa=(X,\Ff,B,\Mm,t)$ and let $g: \widetilde{X}\rightarrow X$ be a foliated log resolution of $\Aa$, such that the induced birational maps $g_i: \widetilde{X}\dashrightarrow X_i$ are morphisms. Let $0<\delta\ll 1$ be a real number such that
    $$K_{g^{-1}_*\Aa}+(1-2\delta)(\Exc(g)^{\ninv}+(1-t)\Exc(g)^{\inv})\geq K_{g^*\Aa}$$
   and let $\widetilde{\Aa}:=(g^{-1}_*\Aa,(1-\delta)(\Exc(g)^{\ninv}+(1-t)\Exc(g)^{\inv})$. Then $\widetilde{\Aa}$ is klt. Let $f: \bar\Aa'\rightarrow\widetilde{\Aa}$ be a $\mathbb Q$-factorial terminalization of $\widetilde{\Aa}$, whose existence is guaranteed by Definition-Theorem \ref{defthm: terminalization}, and let $\bar\Aa:=(\bar\Aa',\epsilon\Exc(f))$ for some $0<\epsilon\ll 1$.
   
   We show that $\bar\Aa$ and $ h=f \circ g $ satisfies our requirements. By Lemma \ref{lem: perturbation terminal}, $\bar\Aa$ is terminal. By our construction, $\phi_i$ is $K_{\Aa}$-negative. Since $h$ is $K_{\bar\Aa}$-negative, $h_i$ is $K_{\bar\Aa}$-negative, and we have
   $$K_{\bar{\Aa}}=h_i^*K_{\Aa_i}+E_i$$
   for some $E_i\geq 0$ that is exceptional$/X_i$ and $\Supp E_i=\Exc(h_i)$. By \cite[Theorem 2.1.1]{CHLMSSX25}, we may run a $K_{\bar\Aa}$-MMP$/X_i$ 
   $$\bar\Aa\dashrightarrow \Aa_{i,1}\dashrightarrow \Aa_{i,2}\dashrightarrow\dots\dashrightarrow \Aa_{i,j}\dashrightarrow\dots$$
   with scaling of ample divisor which either terminates, or the limit of the scaling numbers is $0$. In either case, there exists $j_i$ such that $K_{\Aa_{i,j_i}}$ is a limit of movable$/X_i$ $\mathbb R$-divisors. Since $K_{\Aa_{i,j_i}}\sim_{\mathbb R,X_i}E_{i,j_i}$ where $E_{i,j_i}$ is the image of $E_i$ on $X_{i,j_i}$, by  \cite[Lemma 3.3]{Bir12}, we have $E_{i,j_i}=0$. Thus $X_{i,j_i}=X_i$. Finally, note that
   $$(h_i)_*\bar\Aa=(\phi_i)_*\Aa=\Aa_i.$$
   We conclude the proof.
\end{proof}

\subsection{Finiteness of weak lc models}

\begin{defn}\label{defn: polytope afs}
Let $\pi\colon X\rightarrow U$ be a projective morphism between normal quasi-projective varieties. A \emph{tuple} of nef$/U$ $\bb$-Cartier $\bb$-divisors on $X$ is of the form $\MM:=(\Mm_1,\dots,\Mm_n)$ for some positive integer $n$ where each $\Mm_i$ is a nef$/U$ $\bb$-Cartier $\bb$-divisor. We define $\dim\MM:=n$. We define
$$\Span_{\mathbb R}(\MM):=\left\{\sum_{i=1}^n a_i\Mm_i\middle| a_i\in\mathbb R\right\}, \Span_{\mathbb R_{\geq 0}}(\MM):=\left\{\sum_{i=1}^n a_i\Mm_i\middle| a_i\in\mathbb R_{\geq 0}\right\}.$$
For any nef$/U$ $\bb$-Cartier $\bb$-divisor $\NN$, we define $(\MM,\NN):=(\Mm_1,\dots,\Mm_n,\NN)$.

Given a tuple $\MM$ of nef$/U$ $\bb$-Cartier $\bb$-divisors on $X$, let $V\subset\Weil_{\mathbb R}(X)\times\Span_{\mathbb R}(\MM)$ be a finite dimensional affine subspace, $\Ff$ a foliation on $X$, and $t\in [0,1]$ a real number. We define:
    $$V_{(X,\Ff,t)}:=\{\Aa\mid \Aa=(X,\Ff,B,\Mm,t),(B,\Mm)\in V\}$$
    $$\mathcal{L}(V_{(X,\Ff,t)}):=\{\Aa \mid \Aa=(X,\Ff,B,\Mm,t)\in V_{(X,\Ff,t)},\Aa\text{ is lc},\Mm\in\Span_{\mathbb R_{\geq 0}}(\MM)\}.$$
    $$\mathcal{E}_{\pi}(V_{(X,\Ff,t)}):=\{\Aa\mid\Aa\in \mathcal{L}(V_{(X,\Ff,t)}), K_{\Aa}\text{ is pseudo-effective}/U\}.$$
\end{defn}

\begin{lem}\label{lem: finite non-positive face}
Let $\Aa/U:=(X,\Ff,B,\Mm,t)/U$ be a $\mathbb Q$-factorial klt algebraically integrable adjoint foliated structure such that $B+\Mm_X$ is big$/U$. Then there are finitely many $K_{\Aa}$-non-positive extremal faces$/U$. Moreover, all these faces are contractible.
\end{lem}
\begin{proof}
     By \cite[Lemma 3.26]{CHLMSSX25} there exists an ample $\mathbb R$-divisor $A$ on $X$ and a klt algebraically integrable adjoint foliated structure $\Aa'/U:=(X,\Ff,B',\Mm',t)/U$ such that $K_{\Aa}\sim_{\mathbb R,U}K_{\Aa'}+A$. Thus any $K_{\Aa}$-non-positive extremal face$/U$ is a $(K_{\Aa'}+\frac{1}{2}A)$-negative extremal face$/U$. By \cite[Theorem 1.3]{CHLMSSX24}, there are only finitely many $(K_{\Aa'}+\frac{1}{2}A)$-negative extremal rays$/U$, hence there are only finitely many $(K_{\Aa'}+\frac{1}{2}A)$-negative extremal faces$/U$. 
     
     These faces are rational extremal faces by \cite[Theorem 1.3(4)]{CHLMSSX24}, so each face has a supporting function$/U$ $H$ such that $H-(K_{\Aa'}+\frac{1}{2}A)$ is ample$/U$. The existence of the contraction follows from the base-point-freeness theorem \cite[Theorem 2.1.2]{CHLMSSX25}.
\end{proof}

\begin{thm}[Finiteness of weak lc models]\label{thm: finiteness of weak lc models}
Assume that
\begin{itemize}
    \item $\pi\colon X\rightarrow U$ is a projective morphism between normal quasi-projective varieties,
    \item $\MM$ is a tuple of nef$/U$ $\bb$-Cartier $\bb$-divisors,
    \item $V$ is a finite dimensional rational affine subspace of $\Weil_{\mathbb R}(X)\times\Span_{\mathbb R}(\MM)$,
    \item $\Ff$ is an algebraically integrable foliation on $X$,
    \item $t\in [0,1]$ is a rational number. 
    \item $\mathcal{C}\subset\mathcal{L}(V_{(X,\Ff,t)})$ is a rational polytope, such that for any $\Aa:=(X,\Ff,B,\Mm,t)\in\mathcal{C}$, $\Aa$ is klt and $B+\Mm_X$ is big$/U$.
\end{itemize}
Then there are finitely many birational maps$/U$ $\psi_j: X\dashrightarrow Z_j$, $1\leq j\leq l$ satisfying the following. For any $\Aa\in\mathcal{C}$ and any weak lc model $\Aa_Z/U$ of $\Aa/U$ with induced birational map $\psi\colon X\dashrightarrow Z$, there exists an index $1\leq j\leq l$ such that $\psi_j\circ\psi^{-1}\colon Z\dashrightarrow Z_j$ is an isomorphism.
\end{thm}
\begin{proof}
By \cite[Theorem 2.5.2]{CHLMSSX25}, there are finitely many birational maps$/U$ $\phi_j: X\dashrightarrow Y_j$, $1\leq j\leq k$ satisfying the following. 
\begin{enumerate}
\item For any $\Aa\in\mathcal{C}$ such that $K_{\Aa}$ is pseudo-effective$/U$, there exists an index $1\leq j\leq k$ such that $(\phi_j)_*\Aa/U$ is a $\mathbb Q$-factorial good minimal model of $\Aa/U$.
\item For any $\Aa\in\mathcal{C}$ and any $\mathbb Q$-factorial minimal model $\Aa_Y/U$ of $\Aa/U$ with induced birational map $\phi\colon X\dashrightarrow Y$, there exists an index $1\leq j\leq k$ such that $\phi_j\circ\phi^{-1}\colon Y\dashrightarrow Y_j$ is an isomorphism.
\end{enumerate}
Let
$$\mathcal{C}_j:=\{(\phi_j)_*\Aa\mid K_{(\phi_j)_*\Aa}\text{ is nef}/U, \Aa\in\mathcal{C}\}.$$
Since nef is a close condition and $\mathcal{C}$ is a rational polytope, by \cite[Theorem 2.5.3]{CHLMSSX25}, $\mathcal{C}_j$ is a rational polytope.

Now for any $\Aa\in\mathcal{C}$ and any weak lc model  $\Aa_Z/U$ of $\Aa/U$ with induced birational map $\psi\colon X\dashrightarrow Z$, we let
$$\mathcal{S}:=\{E\mid E\text{ is a prime divisor on }X\text{ that is exceptional}/Z,a(E,\Aa)=a(E,\Aa_Z)\}.$$
By \cite[Theorem 2.2.3]{CHLMSSX25}, there exists a birational morphism $g: Z'\rightarrow Z$ such that $Z'$ is $\Qq$-factorial and the divisors contracted by $g$ are exactly the divisors in $\mathcal{S}$. Let $\Aa_{Z'}:=h^*\Aa_Z$, then $\Aa_{Z'}/U$ is a $\Qq$-factorial good minimal model of $\Aa/U$, hence the induced birational map $\phi: X\dashrightarrow Z'$ is $\phi_j$ for some $1\leq j\leq k$. Now the contraction $g$ is a contraction of a $K_{\Aa_{Z'}}$-trivial extremal face $F$ in $\overline{NE}(Y_j/Z)$ and $\Aa_{Z'}\in\mathcal{C}_j$. By linearity of intersection numbers and since $K_{\Aa_j}$ is nef$/U$ for any $\Aa_j\in\mathcal{C}_j$, there exists a vertex $\Aa_j^0$ of $\mathcal{C}_j$ such that $F$ is a $K_{\Aa_j^0}$-trivial extremal face$/U$. 

There are only finitely many indices $j$. For each $j$, there are only finitely many vertices of $\mathcal{C}_j$. For each vertex $\Aa_j^0$ of $\mathcal{C}_j$, by Lemma \ref{lem: finite non-positive face}, there are only finitely many $K_{\Aa_j^0}$-trivial extremal faces. Thus there are only finitely many possibilities of $Z$.
\end{proof}

\section{Sarkisov program}\label{sec: sarkisov program}

\subsection{Notation and set-up}
\begin{defn}\label{defn: sarkisov links}
Let $f_1: X_1\rightarrow Z_1$ and $f_2: X_2\rightarrow Z_2$ be two $\mathbb Q$-factorial Mori fiber spaces. A \emph{Sarkisov link of Type} I (resp. II, III, IV) between $f_1: X_1\rightarrow Z_1$ and $f_2: X_2\rightarrow Z_2$ is of the following form:
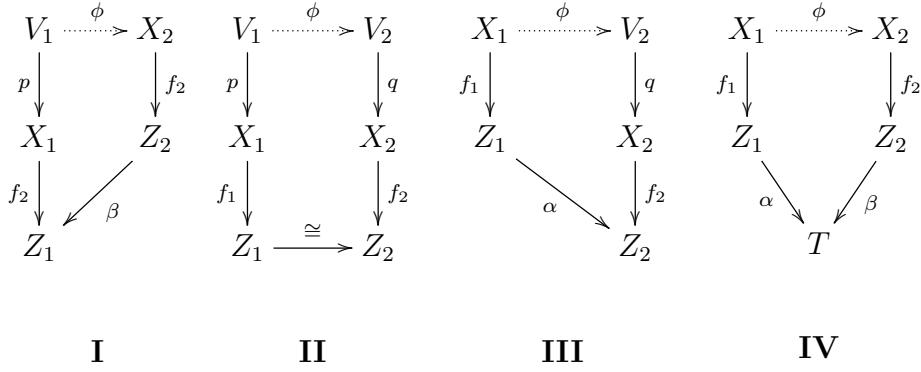
\begin{center}
  $\xymatrix@C=5pt{
    V_1\ar[d]_p\ar@{.>}^{\phi}[rr]&&X_2\ar[d]^{f_2}\\
    X_1\ar[d]_{f_2}& &Z_2\ar[dll]^{\beta}\\
Z_1 &\\
  & \textbf{I} &}$
  $\xymatrix@C=5pt{
    V_1\ar[d]_p\ar@{.>}^{\phi}[rr]& &V_2\ar[d]^{q}&\\
    X_1\ar[d]_{f_1} &&X_2\ar[d]^{f_2}\\
  Z_1\ar[rr]^{\cong}&&Z_2\\
  & \textbf{II} & }$
    $
  \xymatrix@C=5pt{
    X_1\ar@{.>}^{\phi}[rr]\ar[d]_{f_1}&& V_2\ar[d]^q& \\
    Z_1\ar[rrd]_{\alpha}         && X_{2}\ar[d]^{f_{2}}&\\
   & &Z_2\\
    & \textbf{III} &
  }
  $
    $\xymatrix@C=5pt{
    X_1\ar[d]_{f_1}\ar@{.>}^{\phi}[rr]&&X_2\ar[d]^{f_2}\\
    Z_1\ar[dr]_{\alpha}&&Z_2\ar[dl]^{\beta}\\
  &T &\\
  &  \textbf{IV}&}$   
\end{center}
Here $p,q,\alpha,\beta$ are contractions of relative Picard number $1$, $p,q$ are birational, all diagrams commute, and $\phi$ is a sequence of flips. Moreover, if all maps and morphisms in the commutative diagram are projective over a base $U$, then we call the Sarkisov links as Sarkisov links$/U$.
\end{defn}

\begin{setup}\label{setup: setup sarkisov link}
$\Aa_W,W,\Ff_W,B_W,\Mm,t,U,g,h,\Aa,X,\Ff,B,\Aa',X',\Ff',B',f,f',L,R,\Ll,\RR,\Aa_W(\cdot,\cdot)$, $K_W(\cdot,\cdot)$ are defined in the following way.
    \begin{enumerate}
        \item $\Aa_W/U=(W,\Ff_W,B_W,\Mm,t)$ is a $\mathbb Q$-factorial terminal algebraically integrable adjoint foliated structure such that $\Mm$ is NQC$/U$ and descends to $W$,  and $t\in\mathbb Q\cap [0,1)$.
        \item $g: W\rightarrow X$ and $g': W\rightarrow X'$ are projective birational morphisms$/U$.
        \item $\Aa:=g_*\Aa_W=:(X,\Ff,B,\Mm,t)$ and $\Aa':=g'_*\Aa_W=:(X',\Ff',B',\Mm,t)$.
        \item $g$ and $g'$ are $K_{\Aa_W}$-negative and $X,X'$ are $\mathbb Q$-factorial.
        \item $f: X\rightarrow Z$ is a $K_{\Aa}$-Mori fiber space$/U$ and $f': X'\rightarrow Z'$ is a $K_{\Aa'}$-Mori fiber space$/U$.
        \item $L$ (resp. $R$) is an ample$/U$ $\mathbb R$-divisor on $X$ (resp. $X'$) such that $K_{\Aa}+L\sim_{\mathbb R}f^*A_Z$ and  $K_{\Aa'}+R\sim_{\mathbb R}f'^*A_{Z'}$ for some ample$/U$ $\mathbb R$-divisors $A_Z$ and $A_{Z'}$.
        \item $\Ll:=\overline{L}$ and $\RR:=\overline{R}$.
        \item For any non-negative real numbers $a,b$, we denote by $\Aa_W(a,b):=(\Aa_W,a\Ll+b\RR)$ and $K_W(a,b):=K_{\Aa_W(a,b)}$.
    \end{enumerate}
\end{setup}

\begin{Asp}\label{Asp: auxiliary setup}
     Notation and conditions as in Set-up \ref{setup: setup sarkisov link}. $h_i$, $X_i,\Aa_i,\Aa_i(\cdot,\cdot),K_i(\cdot,\cdot),f_i$, $l_i,r_i,Z_i$ are as follows: 
\begin{enumerate}
    \item  $0\leq l_i\leq 1$ and $0\leq r_i\leq 1$ are two real numbers.
    \item  $h_i: W\dashrightarrow X_i$ is a map$/U$.
    \item  $\Aa_i:=(h_i)_*\Aa_W$, $\Aa_i(a,b):=(h_i)_*\Aa_W(a,b)$, and $K_i(a,b):=(h_i)_*K_W(a,b)$ for any $a,b\geq 0$. In particular, $K_i(a,b)=K_{\Aa_i(a,b)}$.
    \item $f_i: X_i\rightarrow Z_i$ is a $\mathbb Q$-factorial $K_{\Aa_i}$-Mori fiber space$/U$.
    \item  $h_i$ is $K_W(l_i,r_i)$-non-positive and $\Aa_i(l_i,r_i)$ is $\mathbb Q$-factorial klt. In particular, $h$ does not extract any divisor.
    \item $K_i(l_i,r_i)$ is nef$/U$ and $K_i(l_i,r_i)\sim_{\mathbb R,Z_i}0$.
\end{enumerate}
\end{Asp}

\begin{defn}\label{defn: setup sarkisov link, induction}
   Notation and conditions as in Set-up \ref{setup: setup sarkisov link}. Assume Assumption \ref{Asp: auxiliary setup} holds for some $i\in\mathbb N$. We define $C_i,\mu_i,\Aa_{W,i}(\cdot),\Aa_i(\cdot),K_{W,i}(\cdot),K_i(\cdot),\Ii_i,s_i,l_i',r_i',\Aa_{W,i},\Aa'_{W,i},\Aa_i^0,\Aa_i'$, $K_{W,i},K_{W,i}',K_i,K_i',\Aa_{W,i}'(\epsilon)$, $\Aa_i'(\epsilon),K_{W,i}'(\epsilon),K_i'(\epsilon)$ in the following way.
\begin{enumerate}
    \item $C_i$ is a general $f_i$-vertical curve.
    \item $\mu_i$ is the unique positive real number such that $\RR_{X_i}\equiv_{Z_i}\mu\Ll_{X_i}$. This is because $\rho(X_i/Z_i)=1$, $\dim X_i>\dim Z_i$, and $\RR_{X_i},\Ll_{X_i}$ are big$/U$.
    \item $\Aa_{W,i}(\lambda):=\Aa_{W}(l_i-\mu_i\lambda,r_i+\lambda)$, $\Aa_i(\lambda):=\Aa_i(l_i-\mu_i\lambda,r_i+\lambda)$, $K_{W,i}(\lambda):=K_{\Aa_{W,i}(\lambda)}$, $K_i(\lambda):=K_{\Aa_i(\lambda)}$.
    \item $\Ii_i$ is the set of all real numbers $\lambda$ satisfying the following.
    \begin{enumerate}
    \item $0\leq \lambda\leq\frac{l_i}{\mu_i}$.
    \item $\Aa_{W,i}(\lambda)\geq \Aa_i(\lambda)$, and
    \item $K_i(\lambda)$ is nef$/U$.
    \end{enumerate}
    \item $s_i:=\sup\{\gamma\mid\gamma\in\Ii_i\}$.
    \item $l_i':=l_i-s_i\mu_i$ and $r_i':=r_i+s_i$.
    \item  $\Aa_{W,i}:=\Aa_{W,i}(0)$, $\Aa_{W,i}':=\Aa_{W,i}(s_i)$, $\Aa^0_{i}:=\Aa_i(0)$, $\Aa_i':=\Aa_i(s_i)$,  $K_{W,i}:=K_{W,i}(0)$, $K_{W,i}':=K_{W,i}(s_i)$, $K_i:=K_i(0)$, $K_i':=K_i(s_i)$.
$\Aa'_{W,i}(\epsilon):=\Aa_{W,i}(s_i+\epsilon)$, $\Aa'_{i}(\epsilon):=\Aa_{i}(s_i+\epsilon)$, $K_{W,i}'(\epsilon):=K_{W,i}(s_i+\epsilon)$, $K_i'(\epsilon):=K_i(s_i+\epsilon)$.
\end{enumerate}
\end{defn}

\begin{lem}\label{lem: sarkisov h<=1}
Notation and conditions as in Set-up \ref{setup: setup sarkisov link}. Assume Assumption \ref{Asp: auxiliary setup} holds for some $i\in\mathbb N$. Adopt notations in Definition \ref{defn: setup sarkisov link, induction}. Then we have:
\begin{enumerate}
    \item Either $\Ii_i=\{0\}$ or $\Ii_i$ is a closed interval, and $s_i\in\Ii$.
    \item $l_i'=l_i$ if and only if $r_i'=r_i$ if and only if $s_i=0$.
    \item $\Ii\subset [0,1-r_i]$. Moreover, $r_i'\leq 1$, and if $r_i'=1$, then $l_i'=0$.
    \item $\Aa_i'$ is klt.
\end{enumerate}
\end{lem}
\begin{proof}
It is obvious that $0\in\Ii$. The values of $\lambda$ which satisfy each of the three conditions of Definition \ref{defn: setup sarkisov link, induction}(3) form a closed and connected set in $\mathbb R_{\geq 0}$, and so is their intersection. This implies (1) and (2) follows from (1).

Assume that (3) does not hold, then $r_i'>1$, or $r_i'=1$ and $l_i'>0$. Then $\Aa_W'$ is terminal as $\RR,\Ll$ descend to $W$. 
By Definition \ref{defn: setup sarkisov link, induction}(4.b) , $\Aa_{W,i}'\geq\Aa_i'$. Thus by Lemma \ref{lem: liu21 3.3}, $\Aa_i'$ is klt. We have
$$0\equiv_{Z_i}K_i'=K_i(0,1)+l_i'\Ll_{X_i}+(r_i'-1)\RR_{X_i}.$$
Since $\Ll_{X_i},\RR_{X_i}$ are big$/U$, $\rho(Y_i/Z_i)=1$, and $\dim X_i>\dim Z_i$, $\Ll_{X_i},\RR_{X_i}$ are ample$/Z_i$. Thus $K_i(0,1)$ is anti-ample$/Z_i$. In particular, $K_i(0,1)$ is not pseudo-effective. This is not possible as 
$$K_W(0,1)=K_{\Aa_{W}}+g'^*R\geq g'^*(K_{\Aa'}+R)\sim_{\mathbb R}(f'\circ g')^*A_{Z'}$$
is pseudo-effective, $K_i(0,1)=(h_i)_*K_W(0,1)$, and $h_i$ does not extract any divisor. This implies (3).

(4) Since $\Aa_{W,i}'\geq\Aa_i'$ and $\Aa_{W,i}'$ is terminal, $\Aa_i'$ is klt.
\end{proof}

\subsection{Construction of Sarkisov link}

\begin{cons}\label{cons: type i and ii}
Notation and conditions as in Set-up \ref{setup: setup sarkisov link}. Assume Assumption \ref{Asp: auxiliary setup} holds for some $i\in\mathbb N$. Adopt notations in Definition \ref{defn: setup sarkisov link, induction}. Assume that $s_i<\frac{l_i}{\mu_i}$ and $\Aa_{W,i}'(\epsilon)\not\geq\Aa_i'(\epsilon)$ for any $0<\epsilon\ll 1$. 

We construct $h_{i+1}$, $X_{i+1},\Aa_{i+1},\Aa_{i+1}(\cdot,\cdot),K_{i+1}(\cdot,\cdot),f_{i+1},l_{i+1},r_{i+1},Z_i$, $p_i,V_i,\phi_i$, and possibly construct $\beta_i,q_i,V_i'$ in the following way, so that Assumption \ref{Asp: auxiliary setup} holds for $i+1$, and $\Aa_i(l_i,r_i)$ and $\Aa_{i+1}(l_i,r_i)$ are crepant.

We define $l_{i+1}:=l_i-s_i\mu_i$ and $r_{i+1}:=r_i+s_i$. 

By our assumption, $K_i'$ is nef$/U$. By Lemma \ref{lem: liu21 3.3}, there exists a prime divisor $E$ on $W$ that is exceptional$/X_i$, such that
$$a(E,\Aa_{W,i}'(\epsilon))>a(E,\Aa_i'(\epsilon))$$
for any $\epsilon>0$. Since $E$ is on $W$, $a(E,\Aa_{W,i}'(\epsilon))\leq 0$ for any $\epsilon>0$, so $a(E,\Aa_i'(\epsilon))<0$ for any $\epsilon>0$, and so
$$a(E,\Aa_i')=\lim_{\epsilon\rightarrow 0^+}a(E,\Aa_i'(\epsilon))\leq 0.$$
By Lemma \ref{lem: sarkisov h<=1}, $\Aa_i'$ is klt. By \cite[Theorem 2.2.3]{CHLMSSX25}, there exists a projective birational morphism $p_i: V_i\rightarrow X_i$ such that $\Exc(p_i)=E$ and $V_i$ is $\mathbb Q$-factorial. Since $X_i$ is $\mathbb Q$-factorial, $\rho(V_i/X_i)=1$. 
Pick $0<\delta\ll\epsilon\ll 1$ and let
$$\widetilde\Aa_i:=\Aa_i(l_{i+1}-\mu_i\epsilon-\delta,r_{i+1}+\epsilon).$$
Since $a(E,\Aa_i'(\epsilon))<0$ and $0<\delta\ll\epsilon$, $a(E,\widetilde\Aa_i)<0$. Let
$$\widetilde\Aa_{V_i}:=p_i^*\widetilde\Aa_i.$$
Then $\widetilde\Aa_{V_i}$ is klt. Moreover, we have
$$K_{\widetilde\Aa_i}\equiv_{Z_i}-\delta\Ll_{X_{i+1}}$$
is anti-big$/Z_i$, hence $K_{\widetilde\Aa_i}$ is not pseudo-effective$/Z_i$, and hence $K_{\widetilde\Aa_{V_i}}$ is not pseudo-effective$/Z_i$. By \cite[Theorem 2.1.1]{CHLMSSX25}, we may run a $K_{\widetilde\Aa_{V_i}}$-MMP$/Z_i$ with scaling of an ample divisor which terminates with a Mori fiber space of $\widetilde\Aa_{V_i}/Z_i$. Let $\Aa_{V_i}:=p_i^*\Aa_i'$, then $K_{\Aa_{V_i}}\sim_{\mathbb R,Z_i}0$, so this MMP is $K_{\Aa_{V_i}}$-trivial. There are two possibilities:

\medskip

\noindent\textbf{Case I.} After finitely many flips $\phi_i: V_i\dashrightarrow X_{i+1}$, we obtain a Mori fiber space$/Z_i$ $f_{i+1}: X_{i+1}\rightarrow Z_{i+1}$. We have $\rho(X_{i+1}/Z_{i+1})=1$ and $\rho(X_{i+1}/Z_i)=2$, so there exists a non-trivial contraction $\beta_i: Z_{i+1}\rightarrow Z_i$. Then we obtain a Sarkisov link$/U$ of Type I:
\begin{center}
  $\xymatrix{
    V_i\ar[d]_p\ar@{.>}^{\phi_i}[rr]&&X_{i+1}\ar[d]^{f_{i+1}}\\
    X_i\ar[d]_{f_i}& &Z_{i+1}\ar[dll]^{\beta_i}\\
Z_i &}$
  \end{center}
 We let $h_{i+1}: W\dashrightarrow X_{i+1}$ be the induced birational map, $\Aa_{i+1}:=(h_{i+1})_*\Aa_W$, $\Aa_{i+1}(a,b):=(h_{i+1})_*\Aa_W(a,b)$, and $K_{i+1}(a,b):=(h_{i+1})_*K_W(a,b)$ for any $a,b\geq 0$. Other notations are defined via Definition \ref{defn: setup sarkisov link, induction} for $i+1$.

\medskip

\noindent\textbf{Case II.} After finitely many flips $\phi_i: V_i\dashrightarrow V_{i}'$, we obtain a divisorial contraction $q_{i}: V_{i}'\rightarrow X_{i+1}$. Then $(q_i\circ\phi_i)_*K_{\widetilde\Aa_{V_i}}$ is not pseudo-effective$/Z_i$, so the induced contraction $f_{i+1}: X_{i+1}\rightarrow Z_{i+1}:\cong Z_i$ is a Mori fiber space.  We let $h_{i+1}: W\dashrightarrow X_{i+1}$ be the induced birational map, $\Aa_{i+1}:=(h_{i+1})_*\Aa_W$, $\Aa_{i+1}(a,b):=(h_{i+1})_*\Aa_W(a,b)$, and $K_{i+1}(a,b):=(h_{i+1})_*K_W(a,b)$ for any $a,b\geq 0$. Other notations are defined via Definition \ref{defn: setup sarkisov link, induction} for $i+1$.

\medskip

We check that Assumption \ref{Asp: auxiliary setup} holds for $i+1$ for both \textbf{Case I} and \textbf{Case II}. By Lemma \ref{lem: sarkisov h<=1}(3), we only need to check Set-up \ref{Asp: auxiliary setup}(4-6). By our construction, $f_{i+1}$ is $K_{i+1}$-trivial, hence $f_{i+1}$ is a $K_{\Aa_{i+1}}$-Mori fiber space as $\Ll_{X_{i+1}},\RR_{X_{i+1}}$ are big. This implies Set-up \ref{Asp: auxiliary setup}(4) for $i+1$.

Since $h_i$ is $K_{W,i}'$-non-positive, $p_i$ is $K_{\Aa_{V_i}}$-trivial, and the induced birational map $W_i\dashrightarrow V_i$ does not extract any divisor, we have that $h_{i+1}$ is $K_{W,i}'$-non-positive. Since $\Aa_i'$ is klt, $\Aa_{V_i}$  and $\Aa_{i+1}^0$ are klt. This implies Set-up \ref{Asp: auxiliary setup}(5) for $i+1$. 

By Lemma \ref{lem: single sarkisov link not nef}, $K_i'\sim_{\mathbb R,Z_i}0$ and $K_i'$ is nef$/U$, so $K_{\Aa_{V_i}}\sim_{\mathbb R,Z_i}0$ and $K_{\Aa_{V_i}}$ is nef$/U$, and so $K_{i+1}\sim_{\mathbb R,Z_i}0$ and $K_{i+1}$ is nef$/U$. This implies Set-up \ref{Asp: auxiliary setup}(6) for $i+1$, and in particular, $\Aa_i'=\Aa_i(l_i,r_i)$, $\Aa_{V_i}$, $\Aa_{i+1}^0=\Aa_{i+1}(l_i,r_i)$ are crepant to each other.
\end{cons}

\begin{lem}\label{lem: no infinite type ii}
Notation and conditions as in Construction \ref{cons: type i and ii}. Assume that we are in {\rm\textbf{Case II}}. Then for any $\delta',\epsilon'>0$, we have the following.
\begin{enumerate}
    \item For any prime divisor $F$ over $X$,
    $$a(F,\Aa_i(l_{i+1}-\mu_i\epsilon'-\delta',r_{i+1}+\epsilon'))\leq a(F,\Aa_{i+1}(l_{i+1}-\mu_i\epsilon'-\delta,r_{i+1}+\epsilon')).$$
    \item There exists a prime divisor $F_0$ over $X$, such that
    $$a(F,\Aa_i(l_{i+1}-\mu_i\epsilon'-\delta',r_{i+1}+\epsilon'))<a(F,\Aa_{i+1}(l_{i+1}-\mu_i\epsilon'-\delta',r_{i+1}+\epsilon')).$$
\end{enumerate}
\end{lem}
\begin{proof}
By Construction \ref{cons: type i and ii}, we know that for some $0<\delta\ll\epsilon<1$, $(q_i\circ \phi_i)$ is a sequence of steps of a $p_i^*K_i(l_{i+1}-\mu_i\epsilon-\delta,r_{i+1}+\epsilon)$-MMP$/Z_i$, and since this MMP contains a divisorial contraction, it is not trivial. We have
$$p_i^*K_i(l_{i+1}-\mu_i\epsilon-\delta,r_{i+1}+\epsilon)\sim_{\mathbb R,Z_i}\frac{\delta}{\delta'}p_i^*K_i(l_{i+1}-\mu_i\epsilon'-\delta',r_{i+1}+\epsilon')$$
for any $\epsilon,\delta>0$. Therefore,  $(q_i\circ \phi_i)$ is a sequence of steps of a $p_i^*K_i(l_{i+1}-\mu_i\epsilon'-\delta',r_{i+1}+\epsilon')$-MMP$/Z_i$, and since this MMP contains a divisorial contraction, it is not trivial. The lemma follows.
\end{proof}

\begin{lem}\label{lem: single sarkisov link not nef}
Notation and conditions as in Set-up \ref{setup: setup sarkisov link}. Assume Assumption \ref{Asp: auxiliary setup} holds for some $i\in\mathbb N$. Adopt notations in Definition \ref{defn: setup sarkisov link, induction}.

Assume that $s_i<\frac{l_i}{\mu_i}$ and $\Aa_{W,i}'(\epsilon)\geq\Aa_i'(\epsilon)$ for any $0<\epsilon\ll 1$. Then $\Aa_i'(\epsilon)$ is klt for any $0\leq\epsilon\ll 1$, and there exists a contraction $\alpha_i: Z_i\rightarrow T_i$ such that $\rho(Z_i/T_i)=1$ and $K_i'(\epsilon)$ is not nef$/T_i$ for any $\epsilon>0$. In particular, $K_i'\sim_{\mathbb R,T_i}0$.
\end{lem}
\begin{proof}
For any $0\leq\epsilon\ll 1$, since $\Aa_{W,i}'(\epsilon)\geq\Aa_i'(\epsilon)$ and $\Aa_{W,i}'(\epsilon)$ is terminal, $\Aa_i'(\epsilon)$ is klt. By our assumption, $K_i'(\epsilon)$ is not nef$/U$ for any $\epsilon>0$.  By Lemma \ref{lem: finite non-positive face}, there are finitely many $K_i'$-trivial extremal rays in $\overline{NE}(X_i/U)$, hence for any $0<\epsilon\ll 1$, there exists a $K_i'(\epsilon)$-negative extremal ray $\Sigma$ that is $K_i'$-trivial. Let $\Sigma'$ be the unique extremal ray of $\overline{NE}(Y_i/Z_i)$, then $\Sigma'$ is $K_i'(\epsilon)$-trivial for any $\epsilon$, so $\Sigma_i'\not=\Sigma_i$, and so $\Sigma_i'$ and $\Sigma_i$ span a $2$-dimensional extremal face $F$ of $\overline{NE}(X_i/U)$. By Lemma \ref{lem: finite non-positive face}, there exists a contraction $X_i\rightarrow T_i$ of $F$ which factors through $Z_i$. We may let $\alpha_i: Z_i\rightarrow T_i$ be the induced contraction.
\end{proof}

\begin{cons}\label{cons: type iii and iv}
Notation and conditions as in Set-up \ref{setup: setup sarkisov link}. Assume Assumption \ref{Asp: auxiliary setup} holds for some $i\in\mathbb N$. Adopt notations in Definition \ref{defn: setup sarkisov link, induction}. Assume that $s_i<\frac{l_i}{\mu_i}$ and $\Aa_{W,i}'(\epsilon)\geq\Aa_i'(\epsilon)$ for any $0<\epsilon\ll 1$. 

We construct $h_{i+1}$, $X_{i+1},\Aa_{i+1},\Aa_{i+1}(\cdot,\cdot),K_{i+1}(\cdot,\cdot),f_{i+1},l_{i+1},r_{i+1},Z_i$, $\alpha_i,T_i,\phi_i$, and possibly construct $\beta_i,q_i,V_i$ in the following way, so that Assumption \ref{Asp: auxiliary setup} holds for $i+1$, and $\Aa_i(l_i,r_i)$ and $\Aa_{i+1}(l_i,r_i)$ are crepant.

We define $l_{i+1}:=l_i-s_i\mu_i$ and $r_{i+1}:=r_i+s_i$.

By Lemma \ref{lem: single sarkisov link not nef}, there exists $\epsilon_0>0$, such that $\Aa_i'(\epsilon)$ is klt for any $0\leq\epsilon\leq\epsilon_0$, and there exists a contraction $\alpha_i: Z_i\rightarrow T_i$ such that $\rho(Z_i/T_i)=1$ and $K_i'(\epsilon)$ is not nef$/T_i$ for any $0<\epsilon\leq\epsilon_0$. We define $\alpha_i,T_i$ as above.

Pick $0<\epsilon\ll\epsilon_0$. Since the generalized boundary of $\Aa_{W,i}'(\epsilon)$ is big$/U$, the generalized boundary of $\Aa_{i}'(\epsilon)$ is big$/U$. By \cite[Theorem 2.1.1]{CHLMSSX25}, we may run a $K_i'(\epsilon)$-MMP$/T_i$ with scaling of an ample divisor which terminates. By Lemma \ref{lem: single sarkisov link not nef}, $K_i'\sim_{\mathbb R,T}0$, so this MMP is $K_i'$-trivial. 

There are three possibilities:

\medskip

\noindent\textbf{Case III.} After finitely many steps of the MMP $\phi_i: X_i\dashrightarrow V_i$ we obtain a divisorial contraction$/T$ $q_i: V_i\rightarrow X_{i+1}$. Then $\rho(X_{i+1}/T_i)=1$ and $(q_i\circ\phi_i)_*K_i'(\epsilon)$ is not pseudo-effective$/T_i$, so the induced contraction $f_{i+1}: X_{i+1}\rightarrow T_i$ is a Mori fiber space. Let $Z_{i+1}:=T_i$, then we obtain a Sarkisov link$/U$ of Type III:
\begin{center}
        $
  \xymatrix{
    X_i\ar@{.>}^{\phi_i}[rr]\ar[d]_{f_i}&& V_i\ar[d]^{q_i}& \\
    Z_i\ar[rrd]_{\alpha_i}         && X_{i+1}\ar[d]^{f_{i+1}}&\\
   & &Z_{i+1}
  }
  $
\end{center}
We let $h_{i+1}: W\dashrightarrow X_{i+1}$ be the induced birational map, $\Aa_{i+1}:=(h_{i+1})_*\Aa_W$, $\Aa_{i+1}(a,b):=(h_{i+1})_*\Aa_W(a,b)$, and $K_{i+1}(a,b):=(h_{i+1})_*K_W(a,b)=K_{\Aa_{i+1}(a,b)}$ for any $a,b\geq 0$.Other notations are defined via Definition \ref{defn: setup sarkisov link, induction} for $i+1$.

\medskip

\noindent\textbf{Case IV.1} After finitely many flips $\phi_i: X_i\dashrightarrow X_{i+1}$ we obtain a Mori fiber space$/T_i$ $f_{i+1}: X_{i+1}\rightarrow Z_{i+1}$. Let $\beta_i: Z_{i+1}\rightarrow T_{i}$ be the induced contraction. Then we obtain a Sarkisov link$/U$ of Type IV:

\begin{center}
       $\xymatrix{
    X_i\ar[d]_{f_i}\ar@{.>}^{\phi_i}[rr]&&X_{i+1}\ar[d]^{f_{i+1}}\\
    Z_i\ar[dr]_{\alpha_i}&&Z_{i+1}\ar[dl]^{\beta_i}\\
  &T_i &}$   
\end{center}
We let $h_{i+1}: W\dashrightarrow X_{i+1}$ be the induced birational map, $\Aa_{i+1}:=(h_{i+1})_*\Aa_W$, $\Aa_{i+1}(a,b):=(h_{i+1})_*\Aa_W(a,b)$, and $K_{i+1}(a,b):=(h_{i+1})_*K_W(a,b)$ for any $a,b\geq 0$. Other notations are defined via Definition \ref{defn: setup sarkisov link, induction} for $i+1$.

\medskip

Before we move on to the last case, we check that for \textbf{Case III} and \textbf{Case IV.1}, Assumption \ref{Asp: auxiliary setup} holds for $i+1$. By Lemma \ref{lem: sarkisov h<=1}(3), we only need to check Set-up \ref{Asp: auxiliary setup}(4-6).

By our construction, $f_{i+1}$ is a $K_{i+1}(\epsilon)$-Mori fiber space, hence $f_{i+1}$ is a $K_{\Aa_{i+1}}$-Mori fiber space as $\Ll_{X_{i+1}},\RR_{X_{i+1}}$ are big. This implies Set-up \ref{Asp: auxiliary setup}(4) for $i+1$.

Since $h_i$ is $K_{W,i}'$-non-positive and $\phi_i$ is $K_i'$-trivial, $h_{i+1}$ is $K_{W,i}'$-non-positive. Since $\Aa_i'$ is klt, so $\Aa_{i+1}^0$ is klt. This implies Set-up \ref{Asp: auxiliary setup}(5) for $i+1$. 

By Lemma \ref{lem: single sarkisov link not nef}, $K_i'\sim_{\mathbb R,T}0$ and $K_i'$ is nef$/U$, so $K_{i+1}\sim_{\mathbb R,T}0$ and $K_{i+1}$ is nef$/U$. This implies Set-up \ref{Asp: auxiliary setup}(6) for $i+1$.

\medskip

\noindent\textbf{Case IV.2.} After finitely many flips $\phi_i: X_i\dashrightarrow X_{i+1}$ we obtain a good minimal model $\Aa_{i+1}(\epsilon)/T_i$ of $\Aa_{i}'(\epsilon)/T_i$. We let $h_{i+1}: W\dashrightarrow X_{i+1}$ be the induced birational map, $\Aa_{i+1}:=(h_{i+1})_*\Aa_W$, $\Aa_{i+1}(a,b):=(h_{i+1})_*\Aa_W(a,b)$, and $K_{i+1}(a,b):=(h_{i+1})_*K_W(a,b)$ for any $a,b\geq 0$. Other notations are defined via Definition \ref{defn: setup sarkisov link, induction} for $i+1$.

Let $C_i'$ be the strict transform of $C_i$ on $X_{i+1}$. By the projection formula, we have
$$K_{i+1}(\epsilon)\cdot C_{i}'=K_{i}'(\epsilon)\cdot C_i=0.$$
We claim that $[C_{i}']$ is an extremal ray in $\overline{NE}(X_{i+1}/T_i)$. Supposes not, then we let $\theta,\tau$ be two different extremal rays of $\overline{NE}(X_{i+1}/T_i)$, and we have $[C_{i}']=a\theta+b\tau$ for some $a,b>0$. Since $K_{i+1}(\epsilon)$ is nef$/T_i$, we have 
$$K_{i+1}(\epsilon)\cdot\theta=K_{i+1}(\epsilon)\cdot\tau=0,$$
which indicates that $K_{i+1}(\epsilon)\sim_{\mathbb R,T_i}0$. However, since $\phi_i$ only consists of flips, $K_{i}'(\epsilon)\sim_{\mathbb R,T_i}0$, which is not possible.

Therefore, $[C_i']$ is an extremal ray in $\overline{NE}(X_{i+1}/T_i)$. By Lemma \ref{lem: single sarkisov link not nef}, $K_i'\sim_{\mathbb R,T}0$ and $K_i'$ is nef$/U$. Since $h_i$ is $K_W'$-non-positive and $\phi_i$ is $K_i'$-trivial, $h_{i+1}$ is $K_W'$-non-positive. By Lemma \ref{lem: single sarkisov link not nef}, $K_i'\sim_{\mathbb R,T}0$ and $K_i'$ is nef$/U$, so $K_{i+1}\sim_{\mathbb R,T}0$ and $K_{i+1}$ is nef$/U$. By Lemma \ref{lem: sarkisov h<=1}, $\Aa_i'$ is klt, so $\Aa_{i+1}^0$ is klt. Since $\RR_{X_{i+1}}$ is big$/U$, by Lemma \ref{lem: finite non-positive face}, there exists a contraction$/T_i$ $f_{i+1}: X_{i+1}\rightarrow Z_{i+1}$ of $[C_i']$. Let $\beta_i: Z_{i+1}\rightarrow T_i$ be the induced contraction, and we obtain a Sarkisov link$/U$ of Type IV:
\begin{center}
       $\xymatrix{
    X_i\ar[d]_{f_i}\ar@{.>}^{\phi_i}[rr]&&X_{i+1}\ar[d]^{f_{i+1}}\\
    Z_i\ar[dr]_{\alpha_i}&&Z_{i+1}\ar[dl]^{\beta_i}\\
  &T_i &}$   
\end{center}

We check that Assumption \ref{Asp: auxiliary setup} holds for $i+1$. By Lemma \ref{lem: sarkisov h<=1}(3) and what we have discussed above, we only need to check Set-up \ref{Asp: auxiliary setup}(4) for $i+1$. Since $C_{i}$ can be any general  vertical$/Z_i$ curve on $X_i$, $f_{i+1}$ is a Mori fiber space. Since $f_{i+1}$ is $K_{i+1}$-trivial and $\Ll_{X_{i+1}}$, $\RR_{X_{i+1}}$ are big$/U$, $f_{i+1}$ is $K_{\Aa_{i+1}}$-negative, hence a $K_{\Aa_{i+1}}$-Mori fiber space$/U$. This concludes Set-up \ref{Asp: auxiliary setup}(4) for $i+1$.

\medskip

In all cases, $\Aa_i(l_i,r_i)$ and $\Aa_{i+1}(l_i,r_i)$ are crepant because $\phi_i$ is $K_i'$-trivial.
\end{cons}

\begin{lem}\label{lem: case iv1}
 Notation and conditions as in Construction \ref{cons: type iii and iv}. Assume that we are in {\rm\textbf{Case IV.1}}. Then $\mu_{i+1}>\mu_i$.
\end{lem}
\begin{proof}
We have 
$$K_{\Aa_{i+1}}+l_{i+1}\Ll_{X_{i+1}}+r_{i+1}\RR_{X_{i+1}}=K_{i+1}(l_{i+1},r_{i+1})=K_{i+1}\sim_{\mathbb R,Z_{i+1}}0$$
and
$$K_{\Aa_{i+1}}+(l_{i+1}+\mu_i\epsilon)\Ll_{X_{i+1}}+(r_{i+1}-\epsilon)\RR_{X_{i+1}}=K_{i+1}(l_{i+1}+\mu_i\epsilon,r_{i+1}-\epsilon)=(\phi_i)_*K_i'(\epsilon)$$
is anti-ample$/Z_{i+1}$. Thus
$$\mu_i\Ll_{X_{i+1}}-\RR_{X_{i+1}}\equiv (\mu_i-\mu_{i+1})\Ll_{X_{i+1}}$$
is anti-ample$/Z_i$, so $\mu_{i+1}>\mu_i$.
\end{proof}

\begin{lem}\label{lem: case iv2}
 Notation and conditions as in Construction \ref{cons: type iii and iv}. Assume that we are in {\rm\textbf{Case IV.2}}. Then:
 \begin{enumerate}
      \item $\phi_i$ is not an isomorphism.
     \item $\mu_{i+1}=\mu_i$.
    \item For any prime divisor $F$ over $X$ and any $0<\epsilon\ll 1$,
    $$a(F,\Aa_i'(\epsilon))\leq a(F,\Aa_{i+1}(\epsilon)).$$
    \item There exists a prime divisor $F_0$ over $X$, such that
    $$a(F_{0},\Aa_i'(\epsilon))<a(F_{0},\Aa_{i+1}(\epsilon))$$
    for any $0<\epsilon\ll 1$.
 \end{enumerate}
\end{lem}
\begin{proof}
Since $K_i'(\epsilon)$ is not nef$/T_i$ but $(\phi_i)_*K_i'(\epsilon)$ is nef$/T_i$, we get (1). 

Let $\psi_i: Y\rightarrow X_i$ and $\psi_{i+1}: Y\rightarrow X_{i+1}$ be common resolution and let $C_Y$ be an irreducible curve on $Y$ such that $\psi_i(C_Y)=C_i$. Then
\begin{align*}
  0&=K_i'(\epsilon)\cdot C_i=\psi_i^*K_i'(\epsilon)\cdot C_Y=(\psi_i')_*\psi_i^*K_i'(\epsilon)\cdot C_i'=(\phi_i)_*K_i'(\epsilon)\cdot C_i'\\
  &=(K_{i+1}(l_i,r_i)-\mu_i\epsilon\Ll_{X_{i+1}}+\epsilon\RR_{X_{i+1}})\cdot C_i'=(-\mu_i\epsilon\Ll_{X_{i+1}}+\epsilon\RR_{X_{i+1}})\cdot C_i')\\
  &=\epsilon(\mu_{i+1}-\mu_i)\Ll_{X_{i+1}}\cdot C_i',
\end{align*}
so $\mu_i=\mu_{i+1}$, which implies (2).

(3)(4) follow as $\phi_i$ is a sequence of steps of a $K_i'(\epsilon)$-MMP$/U$, and by (1), $\phi_i$ is not an isomorphism.
\end{proof}

\subsection{Termination of Sarkisov program}\label{sec: termination of sarkisov}

\begin{lem}\label{lem: li=0 imply ri=1}
Notation and conditions as in Set-ups \ref{setup: setup sarkisov link}. Assume Assumption \ref{Asp: auxiliary setup} holds for some $i$. If $l_i=0$, then $r_i=1$.
\end{lem}
\begin{proof}
Since $K_i(0,r_i)$ is nef$/U$ and $h_i$ is $K_W(0,r_i)$-non-positive, $K_W(0,r_i)$ is pseudo-effective$/U$, hence 
We have $g'_*K_W(0,r_i)\sim_{\mathbb R,Z'}-(1-r_i)R$ is pseudo-effective$/Z'$. Thus $r_i=1$.
\end{proof}

\begin{lem}\label{lem: li=0 imply conclude}
Notation and conditions as in Set-ups \ref{setup: setup sarkisov link}. Assume Assumption \ref{Asp: auxiliary setup} holds for some $i$. If $(l_i,r_i)=(0,1)$, then $f_i=f'$. In particular, $X_i=X'$ and $Z_i=Z'$.
\end{lem}
\begin{proof}
Let $D$ be a sufficiently ample divisor on $X_i$ and $\Dd:=\overline{D}$. Since $\Aa_W$ is terminal and $f'$ is $K_{\Aa_W}$-negative, $\Aa'$ is klt, hence $(\Aa',\RR+\epsilon\Dd)$ is klt for any $0<\epsilon\ll 1$. Since $D$ is big$/U$, $\Dd_{X'}$ is big$/U$, hence $\Dd_{X'}$ is ample$/Z'$, and so $K_{\Aa'}+R+\epsilon\Dd_{X'}$ is ample for any $0<\epsilon\ll 1$. Since $f'$ is $K_{\Aa_W}$-negative, $f'$ is $(K_{\Aa_W}+\RR_W+\epsilon\Dd_W)$-negative for any $0<\epsilon\ll 1$. But $h_i$ is also  $(K_{\Aa_W}+\RR_W+\epsilon\Dd_W)$-negative, hence both $X_i$ and $X'$ are ample models$/U$ of $K_{\Aa_W}+\RR_W+\epsilon\Dd_W$. Thus $h_i=g'$, $X_i\cong X'$, and $K_i(0,1)=K_{\Aa'}+R\sim_{\mathbb R}f'^*A_{Z'}$. Thus $K_i(0,1)$ is a supporting of the unique extremal ray in $\overline{NE}(X'/Z')$. Since $K_i(0,1)\sim_{\mathbb R,Z_i}0$ and $Z_i\not=X_i$, we have $f_i=f'$ and $Z_i=Z'$. 
\end{proof}

\begin{lem}\label{lem: finite output in sarkisov program}
Notation and conditions as in Set-ups \ref{setup: setup sarkisov link}. Assume Assumption \ref{Asp: auxiliary setup} holds for any $i\in\mathbb N$. Then there are only finitely many possibilities of $h_i$.   
\end{lem}
\begin{proof}
    $\Aa_i(l_i,r_i)/U$ is a weak lc model of $\Aa_{W}(l_i,r_i)/U$, so we are done by Theorem \ref{thm: finiteness of weak lc models}.
\end{proof}

\begin{cons}\label{cons: sarkisov program no termination}
Notation and conditions as in Set-up \ref{setup: setup sarkisov link}. Adopt notations in Definition \ref{defn: setup sarkisov link, induction}. 

We construct $h_i$, $X_i,\Aa_i,\Aa_i(\cdot,\cdot),K_i(\cdot,\cdot),f_i$, $l_i,r_i,Z_i$ for any $i\in\mathbb N$ in the following way.

First, we let $h_0:=g$, $X_0:=X$, $\Aa_0:=\Aa$, $\Aa_0(a,b):=g_*\Aa_W(a,b)$, $K_0(a,b):=g_*K_W(a,b)$, $f_0:=f$, $l_0:=1$, $r_0:=0$, $Z_0:=Z$. Assumption \ref{Asp: auxiliary setup} holds for $i=0$.

Then, assume that Assumption \ref{Asp: auxiliary setup} holds for $i$ and $\mu_i<\frac{l_i}{r_i}$. Then we construct $h_{i+1}$, $X_{i+1},\Aa_{i+1},\Aa_{i+1}(\cdot,\cdot),K_{i+1}(\cdot,\cdot),f_{i+1}$, $l_{i+1},r_{i+1},Z_{i+1}$ via Constructions \ref{cons: type i and ii} and \ref{cons: type iii and iv} so that Assumption \ref{Asp: auxiliary setup} holds for $i+1$.

Finally, assume that Assumption \ref{Asp: auxiliary setup} holds for $i=N-1$ and $s_i=\frac{l_i}{\mu_i}$. We let $h_{N}:=h_i$, $X_{N}:=X_i$, $\Aa_{N}:=\Aa_i$, $\Aa_{N}(a,b):=\Aa_i(a,b)$, $K_{N}(a,b):=K_i(a,b)$, $f_{N}:=f_i$, $l_{N}:=0$, $r_{N}:=r_i+s_i$, $Z_{N}:=Z_i$.  Assumption \ref{Asp: auxiliary setup} holds for $N=i+1$ as well.
\end{cons}

\begin{thm}\label{thm: sarkisov terminate}
Notation and conditions as in Construction \ref{cons: sarkisov program no termination}.

Assume that $l_N=0$ for some $N$. Then $f_N=f'$. In particular, $X_N=X'$ and $Z_N=Z'$.
\end{thm}
\begin{proof}
    It follows from Lemmas \ref{lem: li=0 imply ri=1} and \ref{lem: li=0 imply conclude}.
\end{proof}

\begin{lem}\label{lem: l and r stabilize}
    Notation and conditions as in Construction \ref{cons: sarkisov program no termination}. 
    
    Then there exists $i_0\in\mathbb N^+$, such that $l_i=l_{i_0}$ and $r_i=r_{i_0}$ for any $i\geq i_0$.
\end{lem}
\begin{proof}
By Theorem \ref{thm: sarkisov terminate}, we may assume that $l_i>0$ for any $i$. By Lemma \ref{lem: finite output in sarkisov program}, there exists an infinite sequence of positive integers $i_0<i_1<\dots<i_k<\dots$ such that $h_{i_j}=h_{i_{k}}$ for any $j,k$. Thus $l_{i_j}=l_{i_k}$ by the construction of $l_{i_j}$ and $l_{i_k}$ in Constructions \ref{cons: type i and ii} and \ref{cons: type iii and iv}. By Constructions \ref{cons: type i and ii} and \ref{cons: type iii and iv} again, we have that $l_i$ is a decreasing sequence and $r_i$ is an increasing sequence, so $l_i=l_{i_0}$ and $r_i=r_{i_0}$ for any $i\geq i_0$.  
\end{proof}

\begin{thm}\label{thm: termination of sarkisov program}
    Notation and conditions as in Construction \ref{cons: sarkisov program no termination}. Then $l_N=0$ for some $N$.
\end{thm}
\begin{proof}
Suppose that $l_i>0$ for all $i$.
By Lemma \ref{lem: l and r stabilize}, $l:=l_i$ and $r:=r_i$ are constants for $i\gg 0$. By Lemma \ref{lem: finite output in sarkisov program}, there exists an infinite sequence of positive integers $i_0<i_1<\dots<i_k<\dots$ such that $h_{i_j}=h_{i_{k}}$ for any $j,k$.

First we suppose that $\phi_i$ is a Sarkisov link$/U$ of Type I or II for any $i\gg 0$. Since $h_i$ does not extract any divisor, we have $\rho(W/U)\geq\rho(X_i/U)$ for any $i$. Therefore, for $i\gg 0$, $\phi_i$ is a Sarkisov link$/U$ of Type II. By Lemma \ref{lem: no infinite type ii}, there exists $k>j$ and a prime divisor $F_0$ over $X$ such that
$$a(F_0,\Aa_{i_j}(l-\delta,r))<a(F_0,\Aa_{i_k}(l-\delta,r))$$
for any $0<\delta\ll 1$. This is not possible as $\Aa_{i_j}(l-\delta,r)=\Aa_{i_{k}}(l-\delta,r)$.

Therefore, there exists an infinite sequence of positive integers $i_0'<i_1'<\dots<i_k'<\dots$ such that $\phi_{i_j'}$ is a Sarkisov link$/U$ of Type III or IV for any $j$. By Constructions \ref{cons: type i and ii} and \ref{cons: type iii and iv}, we have that $\Aa_i(l,r)$ and $\Aa_{i+1}(l,r)$ are crepant for any $i\gg 0$. Since $\phi_{i_j}$ is a Sarkisov link$/U$ of Type III or IV for any $j$, by Constructions \ref{cons: type i and ii} and \ref{cons: type iii and iv}, $\Aa_{W,i_1}'(\epsilon)\geq\Aa_{i_1}'(\epsilon)$ for any $0<\epsilon\ll 1$, hence $\Aa_{W,i}'(\epsilon)\geq\Aa_{i}'(\epsilon)$ for any $0<\epsilon\ll 1$ and any $i\geq i_1$, and so $\phi_{i}$ is a Sarkisov link$/U$ of Type III or IV for any $i\geq i_0'$. Since $\rho(X_i/U)>0$ for any $i$, $\phi_{i}$ is a Sarkisov link$/U$ of Type IV for any $i\gg 0$.

Since $h_{i_j}=h_{i_{k}}$ for any $j,k$, $\mu_{i_j}=\mu_{i_k}$ for any $j,k\gg 0$. By Lemmas \ref{lem: case iv1} and \ref{lem: case iv2}(2),  $\phi_{i}$ is a Sarkisov link$/U$ of Type IV.2 for any $i\gg 0$. By Lemma \ref{lem: no infinite type ii}, there exists $k>j$ and a prime divisor $F_0$ over $X$ such that
$$a(F_0,\Aa_{i_j}(\epsilon))<a(F_0,\Aa_{i_k}(\epsilon))$$
for any $0<\epsilon\ll 1$. This is not possible.
\end{proof}

\section{Proof of the main theorems}

\begin{proof}[Proof of Theorem \ref{thm: main afs}]
Since $\phi_i$ is also a $((1-s)K_{\Aa_W}+sK_W+H)$-MMP$/U$ for some $0<s\ll 1$ and ample $\mathbb R$-divisor $H$, by \cite[Lemma 3.29]{CHLMSSX25}, possibly replacing $\Aa_W$, we may assume that $\Mm$ is NQC$/U$, $t\in [0,1)$, and $\Aa$ is klt. By Lemma \ref{lem: lift mfs mmp}, we may assume that $\Aa_W$ is $\mathbb Q$-factorial terminal, the nef part of $\Aa_W$ descends to $W$, and $g,g'$ are $K_{\Aa_W}$-negative morphisms. Now we are in the setting of Set-up \ref{setup: setup sarkisov link}. Proceed the construction as in Section \ref{sec: sarkisov program}, we are done by Theorems \ref{thm: termination of sarkisov program} and \ref{thm: sarkisov terminate}.
\end{proof}

\begin{proof}[Proof of Theorem \ref{thm: main algint}]
It is a special case of Theorem \ref{thm: main afs}.
\end{proof}

\begin{proof}[Proof of Theorem \ref{thm: main dim 3 rank 1}]
$K_{\Ff}$ is not pseudo-effective, so it is algebraically integrable by \cite[Theorem 3.1]{LLM23} (essentially \cite[Theorem 1.1]{CP19}). Now Theorem \ref{thm: main dim 3 rank 1} is a special case of Theorem \ref{thm: main algint}.
\end{proof}

\begin{proof}[Proof of Theorem \ref{thm: main dim 3}]
It follows from Theorem \ref{thm: main dim 3 rank 1} and \cite[Theorem 1.1]{Mas24}. (Note that \cite[Theorem 1.1]{Mas24} requires that the foliation is $\mathbb Q$-factorial F-dlt rather than only assuming that the ambient variety is $\mathbb Q$-factorial klt. On the other hand, any output of an MMP of a rank $2$ foliated log smooth foliation on a threefold is $\mathbb Q$-factorial F-dlt.)
\end{proof}


\begin{thebibliography}{99}
\bibitem[ACSS21]{ACSS21} F. Ambro, P. Cascini, V. V. Shokurov, and C. Spicer, \textit{Positivity of the moduli part}, arXiv:2111.00423.

\bibitem[BFSZ24]{BFSZ24} F. Bernasconi, A. Fanelli, J. Schneider, and S. Zimmermann, \textit{Explicit Sarkisov program for regular surfaces over arbitrary fields and applications}, arXiv:2404.03281.

\bibitem[BLZ21]{BLZ21} J. Blanc, S. Lamy, and S. Zimmermann, \textit{Quotients of higher dimensional Cremona groups}, Acta Math. \textbf{226} (2021), no. 2, 211--318.

\bibitem[Bir12]{Bir12} C. Birkar, \textit{Existence of log canonical flips and a special LMMP}, Publ. Math. IHÉS \textbf{115} (2012), 325--368.

\bibitem[BCHM10]{BCHM10} C. Birkar, P. Cascini, C. D. Hacon, and J. M\textsuperscript{c}Kernan, \textit{Existence of minimal models for varieties of log general type}, J. Amer. Math. Soc. \textbf{23} (2010), no. 2, 405--468.


\bibitem[BZ16]{BZ16} C. Birkar and D.-Q. Zhang, \textit{Effectivity of Iitaka fibrations and pluricanonical systems of polarized pairs}, Publ. Math. IHÉS \textbf{123} (2016), 283--331.


\bibitem[Bru15]{Bru15} M. Brunella, \textit{Birational geometry of foliations}, IMPA Monographs \textbf{1} (2015), Springer, Cham.

\bibitem[BM97]{BM97} A. Bruno and K. Matsuki, \textit{Log Sarkisov proram}, Int. J.Math. \textbf{8} (1997), no. 4, 451--494.

\bibitem[CP19]{CP19} F. Campana and M. P\u{a}un, \textit{Foliations with positive slopes and birational stability of orbifold cotangent bundles}, Publ. Math. IHÉS \textbf{129} (2019), 1--49.

\bibitem[CHLMSSX24]{CHLMSSX24} P. Cascini, J. Han, J. Liu, F. Meng, C. Spicer, R. Svaldi, and L. Xie, \textit{Minimal model program for algebraically integrable adjoint foliated structures}, arXiv:2408.14258.

\bibitem[CHLMSSX25]{CHLMSSX25} P. Cascini, J. Han, J. Liu, F. Meng, C. Spicer, R. Svaldi, and L. Xie, \textit{On finite generation and boundedness of adjoint foliated structures}, arXiv:2504.10737.

\bibitem[CS20]{CS20} P. Cascini and C. Spicer, \textit{On the MMP for rank one foliations on threefolds}, arXiv:2012.11433.

\bibitem[CS21]{CS21} P.~Cascini and C. Spicer, \textit{MMP for co-rank one foliations on threefolds}, Invent. Math. \textbf{225} (2021), no. 2, 603--690.
\bibitem[CHLX23]{CHLX23} G. Chen, J. Han, J. Liu, and L. Xie, \textit{Minimal model program for algebraically integrable foliations and generalized pairs}, arXiv:2309.15823. 

\bibitem[CW25]{CW25} Y. Chen and Y. Wang, \textit{A Note on the Sarkisov Program}, In: \textit{Higher Dimensional Algebraic Geometry: A Volume in Honor of V. V. Shokurov} (C. D. Hacon and C. Xu eds), London Mathematical Society Lecture Note Series. Cambridge University Press (2025), 231--263.

\bibitem[Cor95]{Cor95} A. Corti, \textit{Factoring birational maps of threefolds after Sarkisov}, J. Algebraic Geom., \textbf{4} (1995),  no. 2, 223--254.

\bibitem[Dru21]{Dru21} S. Druel, \textit{Codimension 1 foliations with numerically trivial canonical class on singular spaces}, Duke Math. J. \textbf{170} (2021), no. 1, 95--203.

\bibitem[Flo20]{Flo20} E. Floris, \textit{A note on the G-Sarkisov program}, Enseign. Math. \textbf{66} (2020), no. 1--2, 83--92.

\bibitem[FP22]{FP22} E. Floris and B. Pasquier, \textit{A description of the sarkisov program of horospherical varieties via moment polytopes}, arXiv:2212.10304.

\bibitem[Hac12]{Hac12} C. D. Hacon; \textit{The {{Minimal}} model program for {{varieties}} of log general type}. On the webpage of Hacon, \url{https://www.math.utah.edu/~hacon/MMP.pdf}

\bibitem[HL23]{HL23} C. D. Hacon and J. Liu, \textit{Existence of flips for generalized lc pairs}, Camb. J. Math. \textbf{11} (2023), no. 4, 795--828.  

\bibitem[HM13]{HM13} C. D. Hacon and J. M\textsuperscript{c}Kernan, \textit{The Sarkisov program}, J. Algebraic Geom. \textbf{22} (2013), 389--405.

\bibitem[HL22]{HL22} J. Han and Z. Li, \textit{Weak Zariski decompositions and log terminal models for generalized polarized pairs}, Math. Z. \textbf{302} (2022), 707--741.

\bibitem[He24]{He24} Y. He, \textit{On the strong Sarkisov program}, arXiv:2311.08750.

\bibitem[Isk91]{Isk91} V. A. Iskovskikh, \textit{Generators in the two-dimensional Cremona group over a nonclosed field}, Translation of the 1991 paper from Trudy Mat. Inst. Steklov (1991), 173--188.

\bibitem[Isk96]{Isk96} V. A. Iskovskikh, \textit{Factorization of birational mappings of rational surfaces from the point of view of Mori theory}, Uspekhi Mat. Nauk \textbf{51} (1996), no. 4(310), 3--72.

\bibitem[IKT93]{IKT93} V. A. Iskovskikh, F. K. Kabdykairov, and S. L. Tregub, \textit{Relations in a two-dimensional Cremona
group over a perfect field}, Izv. Ross. Akad. Nauk Ser. Mat., \textbf{57} (1993), no. 3, 3--69. Translation in
Russian Acad. Sci. Izv. Math. \textbf{42} (1994), no. 3, 427--478.

\bibitem[Kal13]{Kal13} A.-S. Kaloghiros, \textit{Relations in the Sarkisov program}, Compos. Math. \textbf{149} (2013), no. 10, 1685--1709.

\bibitem[KM98]{KM98} J. Koll\'{a}r and S. Mori, \textit{Birational geometry of algebraic varieties}, Cambridge Tracts in Math. \textbf{134} (1998), Cambridge Univ. Press.

\bibitem[LZ20]{LZ20} S. Lamy and S. Zimmermann, \textit{Signature morphisms from the Cremona group over a non-closed field}, J. Eur. Math. Soc. \textbf{22} (2020), no. 10, 3133--3173.

\bibitem[Liu21]{Liu21} J. Liu, \textit{Sarkisov program for generalized pairs}, Osaka J. Math. \textbf{58} (2021), 899--920.

\bibitem[LLM23]{LLM23} J. Liu, Y. Luo, and F. Meng, \textit{On global ACC for foliated threefolds},  Trans. Amer. Math. Soc. \textbf{376} (2023), no. 12, 8939--8972.

\bibitem[LMX24]{LMX24} J. Liu, F. Meng, and L. Xie, \textit{Minimal model program for algebraically integrable foliations on klt varieties}, arXiv:2404.01559.

\bibitem[Mas24]{Mas24} R. Mascharak, \textit{On the log Sarkisov program for foliations on projective 3-folds}, arXiv:2406.09434.

\bibitem[McQ08]{McQ08} M. McQuillan, \textit{Canonical models of foliations}, Pure Appl. Math. Q. \textbf{4} (2008), no. 3, Special Issue: In honor of Fedor Bogomolov, Part 2, 877--1012.

\bibitem[Miy19]{Miy19} K. Miyamoto, \textit{The Sarkisov program on log surfaces}, arXiv:1910.07025.

\bibitem[Sar80]{Sar80} V. G. Sarkisov, \textit{Birational automorphisms of conic bundles}, Izv. Akad. Nauk SSSR Ser. Math., \textbf{44} (1980), no. 4, 918--945, 974.

\bibitem[Sar82]{Sar82} V. G. Sarkisov, \textit{On conic bundle structures}, Izv. Akad. Nauk SSSR Ser. Math. \textbf{46} (1982), no. 2, 371--408, 432.

\bibitem[Spi20]{Spi20} C. Spicer, \textit{Higher dimensional foliated Mori theory}, Compos. Math. \textbf{156} (2020), no. 1, 1--38.

\bibitem[SS22]{SS22} C. Spicer and R. Svaldi, \textit{Local and global applications of the Minimal Model Program for co-rank 1 foliations on threefolds}, J. Eur. Math. Soc. \textbf{24} (2022), no. 11, 3969--4025.

\bibitem[Sti21]{Sti21} L. Stigant, \textit{Mori Fibrations in Mixed Characteristic}, arXiv:2110.06067.

\end{thebibliography}
\end{document}